\documentclass{article}
\usepackage{amsmath}
\usepackage{amsmath,amssymb,latexsym,color}
\usepackage[mathscr]{eucal}
\newtheorem{thm}{Theorem}[section]

\newtheorem{lemma}[thm]{Lemma}
\newtheorem{cor}[thm]{Corollary}

\newtheorem{example}{Example}[section]
\newtheorem{defin}{Definition}[section]

\newtheorem{remark}{Remark}[section]

\newcommand{\proof}{{\it Proof.\quad}}
\newcommand{\qed}{\hfill\Box\medskip}

\usepackage{CJK}
\begin{document}
\begin{CJK*}{GBK}{song}

\renewcommand{\baselinestretch}{1.2}

\title{\bf Cameron-Liebler sets for maximal totally isotropic flats in classical affine spaces}

\author{
Jun Guo\thanks{Corresponding author. E-mail address: guojun$_-$lf@163.com}
\quad Lingyu Wan\thanks{E-mail address: 18631690160@163.com}\\
{\footnotesize Department of Mathematics, Langfang Normal University, Langfang  065000,  China} }
 \date{ }
 \maketitle

\begin{abstract}
Let $ACG(2\nu,\mathbb{F}_q)$ be the $2\nu$-dimensional classical affine
space with parameter $e$ over a $q$-element finite field $\mathbb{F}_q$,
and ${\cal O}_{\nu}$ be the set of all maximal totally isotropic flats in $ACG(2\nu,\mathbb{F}_q)$.
In this paper,  we discuss Cameron-Liebler sets in ${\cal O}_{\nu}$, obtain several equivalent definitions and present some classification results.

\medskip
\noindent {\em AMS classification}: 05B25, 51E20, 05E30, 51E14

 \noindent {\em Key words}: Cameron-Liebler set, Classical affine space, Maximal totally isotropic flat

\end{abstract}

\section{Introduction}
\underline{}Cameron-Liebler sets of lines were first introduced by Cameron and Liebler \cite{Cameron} in their study of collineation groups of PG$(3,q)$. There have been many results for Cameron-Liebler sets of lines in the projective space PG$(3,q)$. See \cite{Metsch2,Metsch3,Rodgers} for classification results,  and \cite{Bruen,De Beule,Feng,Feng2,Gavrilyuk} for the constructions of non-trivial examples. Over the years, there have been many interesting extensions of this result.
See \cite{Blokhuis,De Beule2,Gavrilyuk2,Guo2,Metsch,Rodgers2} for Cameron-Liebler sets of $k$-spaces in PG$(n,q)$,
 \cite{De Boeck2} for Cameron-Liebler classes in finite sets, and  \cite{Filmus} for Cameron-Liebler sets in several classical distance-regular graphs.

One of the main reasons for studying Cameron-Liebler sets is that there are many connections to other geometric and combinatorial objects,
such as blocking sets, intersecting families, linear codes, and association schemes.
Similar problems have been investigated by various researchers under different names: Boolean degree one functions, completely regular
codes of strength $0$ and covering radius $1$, and tight sets, see \cite{Filmus} for more details
on these connections. Recently, De Boeck et al.  \cite{De Boeck,De Boeck0} studied Cameron-Liebler sets of generators in polar spaces, and D'haeseleer et al. \cite{D'haeseleer2,D'haeseleer} studied Cameron-Liebler sets in $AG(n,q)$. Their researches stimulate us to consider Cameron-Liebler sets in classical affine spaces.

Let $\mathbb{F}_q^{2\nu}$ be the $2\nu$-dimensional classical
space with parameter $e$,  and $ACG(2\nu,\mathbb{F}_q)$ be the corresponding $2\nu$-dimensional classical affine
space over a $q$-element finite field $\mathbb{F}_q$, where $e=1,1/2$ or $0$ according to the symplectic, unitary or orthogonal case, respectively, see Section~2.
Let ${\cal O}_{\nu}$ be the set of all maximal totally isotropic flats in $ACG(2\nu,\mathbb{F}_q)$. Then ${\cal O}_{\nu}$ has a structure of
an association scheme ${\mathfrak X}({\cal O}_{\nu})$, which is called the {\it maximal totally isotropic flat scheme} based on ${\cal O}_{\nu}$, see Section~3.
For any vector $\alpha$  in which its positions correspond
to elements in a set, we denote its value on the position corresponding to
an element $a$ by $(\alpha)_a$.
The {\it characteristic vector} $\chi_{{\cal S}}$ of a subset ${\cal S}$ of ${\cal O}_{\nu}$ is the column vector in which its positions correspond to the elements of ${\cal O}_{\nu}$, such that $(\chi_{{\cal S}})_{F}=1$ if $F\in{\cal S}$ and 0 otherwise.
Let $M$ be the incidence matrix with rows indexed with vectors in $\mathbb{F}_q^{2\nu}$
and columns indexed with flats in ${\cal O}_{\nu}$ such that entry $M_{x,F}=1$ if
and only if $x\in F$.

\begin{defin} \label{defin1.1}\rm
 A subset ${\cal L}$ of ${\cal O}_{\nu}$ is called a {\it Cameron-Liebler set} with parameter $x=|{\cal L}|/\prod_{i=1}^{\nu}(q^{i+e-1}+1)$
in ${\cal O}_{\nu}$ if $\chi_{{\cal L}}\in {\rm Im}(M^{\top})$, where $M^{\top}$ is the transpose matrix of $M$.
\end{defin}

In this paper, we consider Cameron-Liebler sets in ${\cal O}_{\nu}$.
The rest of this paper is structured as follows.
In Section~2, we introduce classical affine  spaces and show that the rank of $M$ is $(q^{\nu}-1)(q^{\nu-1+e}+1)+1$ over the real field $\mathbb{R}$.
In Section~3, we discuss the maximal totally isotropic flat  scheme ${\mathfrak X}({\cal O}_{\nu})$ in $ACG(2\nu,\mathbb{F}_q)$  and obtain some useful results for later reference.
In Section~4, we give several equivalent
definitions for Cameron-Liebler sets in ${\cal O}_{\nu}$. In Section~5, we present some classification results for Cameron-Liebler sets in ${\cal O}_{\nu}$.

\section{Classical  affine spaces}

Let $\mathbb{F}_q^{2\nu}$ be the vector space over a finite field $\mathbb{F}_q$ with $q$ elements,
and $P$ be an $m$-dimensional subspace of $\mathbb{F}_q^{2\nu}$. Denote also by $P$ a  matrix of
rank $m$ in which rows span the subspace $P$ and call the matrix $P$ a {\it matrix representation}
of the subspace $P$.

Let
$$K=\left(
    \begin{array}{cc}
      0&I^{(\nu)}\\
      -I^{(\nu)}&0
    \end{array}
  \right),
  $$
where $I^{(\nu)}$ is the identity matrix of order $\nu$.
The {\it symplectic group} of degree $2\nu$ with respect to $K$ over
$\mathbb{F}_q,$ denoted by $Sp_{2\nu}(\mathbb{F}_q):=Sp_{2\nu}(\mathbb{F}_q,K)$, consists of
all $2\nu\times 2\nu$ nonsingular matrices $T$ over $\mathbb{F}_q$
satisfying $TKT^{\top}=K$.
There is an action of $Sp_{2\nu}(\mathbb{F}_q)$ on $\mathbb{F}_q^{2\nu}$
as follows:
$$\begin{array}{ccc}
\mathbb{F}_q^{2\nu}\times Sp_{2\nu}(\mathbb{F}_q) &\rightarrow& \mathbb{F}_q^{2\nu}\\
(x,T) &\mapsto& xT.
\end{array}
$$
The above action induces an action on the set of subspaces of $\mathbb{F}_q^{2\nu}$, i.e. a subspace $P$
is carried by $T\in Sp_{2\nu}(\mathbb{F}_q)$ to the subspace $PT$.
The vector space $\mathbb{F}_q^{2\nu}$
together with the right multiplication action  of
$Sp_{2\nu}(\mathbb{F}_q)$ is called the $2\nu$-dimensional {\it
symplectic space} over $\mathbb{F}_q$ and also denoted by $\mathbb{F}_q^{2\nu}$.
An $m$-dimensional subspace
$P$ in the $2\nu$-dimensional symplectic space $\mathbb{F}_q^{2\nu}$ is said to be of {\it
type} $(m,2s),$ if $PKP^{\top}$ is of rank $2s$. In particular, subspaces
of type $(m,0)$ are called $m$-dimensional {\it totally isotropic
subspaces}, and $\nu$-dimensional totally isotropic subspaces
are called  {\it maximal totally isotropic subspaces}.

Let  $q=q_0^2$, where $q_0$ is a prime power. Then $\mathbb{F}_{q}$ has an {\it involutive automorphism}
$a\mapsto \overline{a}=a^{q_0}$ and the fixed field of this automorphism is $\mathbb{F}_{q_0}$.
Let $$H=\left(
    \begin{array}{cc}
      0&I^{(\nu)}\\
      I^{(\nu)}&  0
    \end{array}
  \right).$$
The {\it unitary group} of degree $2\nu$ with respect to $H$ over $\mathbb{F}_{q},$
denoted by $U_{2\nu}(\mathbb{F}_{q}):=U_{2\nu}(\mathbb{F}_{q},H)$, consists of all $2\nu\times
2\nu$ nonsingular matrices $T$ over $\mathbb{F}_{q}$ satisfying
$TH\overline{T}^{\top}=H$. The vector space $\mathbb{F}_{q}^{2\nu}$
together with the right multiplication action of
$U_{2\nu}(\mathbb{F}_{q})$ is called the $2\nu$-dimensional {\it
unitary space} over $\mathbb{F}_{q}$ and also denoted by $\mathbb{F}_{q}^{2\nu}$. An $m$-dimensional subspace
$P$ in the $2\nu$-dimensional unitary space $\mathbb{F}_{q}^{2\nu}$ is said to be of {\it
type} $(m,r),$ if $PH\overline{P}^{\top}$ is of rank $r$. In particular,
subspaces of type $(m,0)$ are called $m$-dimensional {\it totally isotropic
subspaces}, and $\nu$-dimensional totally isotropic subspaces are called  {\it maximal totally
isotropic subspaces}.

Let $\mathbb{F}_q$ be a finite field of odd characteristic. Let
$$
S_{2\nu}= \left(\begin{array} {cc}
  0&I^{(\nu)} \\
  I^{(\nu)}&0
 \end{array}\right).$$
The {\it orthogonal
group} of degree $2\nu$ with respect to
$S_{2\nu}$ over $\mathbb{F}_q$, denoted by $O_{2\nu}(\mathbb{F}_q):=O_{2\nu}(\mathbb{F}_q,S_{2\nu})$,
consists of all $2\nu\times 2\nu$ nonsingular matrices $T$ over
$\mathbb{F}_q$ satisfying
$TS_{2\nu}T^{\top}=S_{2\nu}$.
 The vector space
$\mathbb{F}_{q}^{2\nu}$ together with the right multiplication
action of $O_{2\nu}(\mathbb{F}_q)$ is called the $2\nu$-dimensional
{\it orthogonal space} over $\mathbb{F}_{q}$ and also denoted by $\mathbb{F}_{q}^{2\nu}$. An $m$-dimensional
subspace $P$ in the
  $2\nu$-dimensional orthogonal space $\mathbb{F}_{q}^{2\nu}$ is said to be of {\it type}
$(m,2s)$ if $PS_{2\nu}P^{\top}$ is
cogredient to
$\mbox{diag}(S_{2s},0^{(m-2s)})$. In particular, subspaces
of type $(m,0)$ are called $m$-dimensional {\it totally isotropic
subspaces}, and $\nu$-dimensional totally isotropic subspaces are called  {\it maximal
totally isotropic subspaces}.

For the rest of this paper, we adopt the following notation.
Let $\mathbb{F}_q^{2\nu}$ be the $2\nu$-dimensional classical
space with parameter $e$, and  $G_{2\nu}$ be the corresponding classical group of
degree $2\nu$, where $e=1,1/2$ or $0$ according to the symplectic, unitary or orthogonal case, respectively.
By \cite[Theorems~3.7,~5.8 and~6.4]{Wan}, each
group $G_{2\nu}$ is transitive on the set of all subspaces of the
same type in $\mathbb{F}_q^{2\nu}$. Suppose that $P$ is a subspace of
type $\omega$ in $\mathbb{F}_q^{2\nu}$. The cosets of
$\mathbb{F}_q^{2\nu}$ relative to $P$ are called
$\omega$-{\it flats}. In particular, $\{x\}$ is a $(0,0)$-flat for every $x\in\mathbb{F}_q^{2\nu}$,
and $(\nu,0)$-flats are called {\it maximal totally isotropic flats}.
The dimension of an $\omega$-flat $P+x$
is defined to be the dimension of the subspace $P$, denoted by
$\dim(P+x)$. An  $\omega_1$-flat $F_1$  is said to be {\it incident} with an
$\omega_2$-flat $F_2$, if $F_1$ contains or is contained in
$F_2$. The vector set $\mathbb{F}_q^{2\nu}$ with all the
$\omega$-flats   and the incidence relation among them defined
above is said to be the $2\nu$-dimensional {\it classical affine
space} with parameter $e$, denoted by $ACG(2\nu,\mathbb{F}_q)$, where $e=1,1/2$ or $0$ according to the symplectic, unitary or orthogonal case, respectively.

\medskip
\noindent{\bf Remarks.}
(1)  Let $PG(2\nu,q)$ be the $2\nu$-dimensional projective space over $\mathbb{F}_q$, and
$H_{\infty}=(0^{(2\nu,1)},I^{(2\nu)})$ be the hyperplane at infinity.
Then we can consider all the subspaces $P$ in $PG(2\nu,q)$ such that $P$ has a matrix representation
 $$P=\left(\begin{array}{cc}
 1 & P_1\\
 0 & P_2
 \end{array}\right),$$
 where $P_1$ is a $1\times2\nu$ matrix and $P_2$ is a  $\dim P$-dimensional subspace 
 of some type in the classical space $\mathbb{F}_q^{2\nu}$ with parameter $e$.
 In this way we obtain  $ACG(2\nu,\mathbb{F}_q)$.

(2) Let $P+x$ be a  maximal totally isotropic flat in $ACG(2\nu,\mathbb{F}_q)$, where $P$ is a maximal totally isotropic space in $\mathbb{F}_q^{2\nu}$ and $x\in\mathbb{F}_q^{2\nu}$. Then the vectors in $P+x$ may not be isotropic. Let
$e_{j}\,(1\leq j\leq2\nu)$ be the row vector in $\mathbb{F}_q^{2\nu}$ in which its
$j$th coordinate is 1 and all other coordinates are 0.
For the unitary case, by \cite[Lemma~5.1]{Wan}, the equation $a+\overline{a}=0$ has exactly $q^{1/2}$ solutions in $\mathbb{F}_q$. Choose
$P=(I^{(\nu)},0^{(\nu)})$ and $x=ae_{\nu+1}.$
If $a+\overline{a}\not=0$, then $x+e_1\in P+x$ but $x+e_1$ is not isotropic.

\medskip
Denote by $F_1\cap F_2$ the
intersection of the flats $F_1$ and $F_2$, and by $F_1\vee F_2$ the
minimum flat containing both $F_1$ and $F_2$.

\begin{lemma}\label{lemma2.1}{\rm(See \cite{Gruenberg,Wan}).}
Let $F_1=V_1+x_1$ and $F_2=V_2+x_2$ be any two flats in $ACG(2\nu,\mathbb{F}_q)$,
where $V_1$ and $V_2$ are two subspaces of $\mathbb{F}_q^{2\nu}$, and $x_1,x_2\in\mathbb{F}_q^{2\nu}$. Then the following hold.
\begin{itemize}
\item[\rm(i)]
$F_1\cap F_2\neq\emptyset$ if and only if $x_2-x_1\in V_1+V_2$. Moreover,
if $F_1\cap F_2\neq\emptyset$, then $F_1\cap
F_2=V_1\cap V_2+x$, where $x\in F_1\cap F_2$.

\item[\rm(ii)]
$F_1\vee F_2=V_1+V_2+\langle x_2-x_1\rangle+x_1$. Moreover,
$\dim(F_1\vee F_2)=\dim F_1+\dim F_2-\dim(V_1\cap V_2)+\varepsilon$,
where $\varepsilon=0$ or 1 according to  $F_1\cap F_2\not=\emptyset$ or $F_1\cap F_2=\emptyset$, respectively.
\end{itemize}
\end{lemma}

The set of all matrices
$$\left( \begin{array} {cc}
  T&0\\
  v&1
 \end{array} \right),
$$
where $T\in G_{2\nu}$ and $v\in \mathbb{F}_q^{2\nu}$, forms a group
under matrix multiplication. This group is said to be the
\emph{classical affine group} of $ACG(2\nu,\mathbb{F}_q)$, denoted
by $AG_{2\nu}$. Define the action of $AG_{2\nu}$ on
$ACG(2\nu,\mathbb{F}_q)$ as follows:
\begin{eqnarray*}
ACG(2\nu,\mathbb{F}_q)\times AG_{2\nu}
&\rightarrow& ACG(2\nu,\mathbb{F}_q)\\
\left(x,\left( \begin{array} {cc}
  T&0\\
  v&1
 \end{array} \right)\right)&\mapsto& xT+v.
\end{eqnarray*}
The above action induces an action on the set of flats in $\mathbb{F}_q^{2\nu}$, i.e. a flat $P+x$
is carried by
$$\left( \begin{array} {cc}
  T&0\\
  v&1
 \end{array} \right)\in AG_{2\nu}$$ to the flat $PT+(xT+v)$.
Let ${\cal O}(\omega;2\nu)$ be the set of all $\omega$-flats in $ACG(2\nu,\mathbb{F}_q)$ for a
given $\omega$. By \cite{WG}, ${\cal O}(\omega;2\nu)$ forms an orbit under the action of $AG_{2\nu}$.
For an extensive and detailed introduction about classical affine spaces, we refer to \cite{Wan,WG}.

Define a graph ${\cal G}$ in which its vertex set is the set of all vectors in $\mathbb{F}_q^{2\nu}$, and two vertices $x$ and $y$ are adjacent if $\{x\}\vee\{y\}$ is a $(1,0)$-flat in $ACG(2\nu,\mathbb{F}_q)$.

\begin{lemma}\label{lemma2.3}{\rm (See \cite{Brouwer2,Brouwer3}).}
If $\nu\geq 1$, then the following hold.
\begin{itemize}
\item[\rm(i)]
For the symplectic case, the graph ${\cal G}$ is a complete graph with $q^{2\nu}$ vertices.
Moreover, the distinct eigenvalues of ${\cal G}$ are $q^{2\nu}-1$ and $-1$,
and the corresponding multiplicities are $1$ and
$q^{2\nu}-1$, respectively.

\item[\rm(ii)]
For the unitary and orthogonal cases, the graph ${\cal G}$ is a strongly regular graph with parameters
$$(q^{2\nu},(q^\nu-1)(q^{\nu+e-1}+1),q^{2(\nu+e-1)}+q^{\nu}-q^{\nu+e-1}-2,q^{\nu+e-1}(q^{\nu+e-1}+1)).$$
Moreover, the distinct eigenvalues of ${\cal G}$ are $(q^{\nu}-1)(q^{\nu+e-1}+1),q^{\nu}-q^{\nu+e-1}-1$ and $-(q^{\nu+e-1}+1)$,
and the corresponding multiplicities are $1,(q^{\nu}-1)(q^{\nu+e-1}+1)$ and
$(q^{\nu}-1)(q^{\nu}-q^{\nu+e-1})$, respectively.
\end{itemize}
\end{lemma}

Let $ACG(2\nu,\mathbb{F}_q)$ be the $2\nu$-dimensional classical affine
space with parameter $e$, where $e=1,1/2$ or $0$ according to the symplectic, unitary or orthogonal case, respectively. Recall that ${\cal O}((m,0);2\nu)$ is the set of all $(m,0)$-flats in $ACG(2\nu,\mathbb{F}_q)$.
For convenience, we write ${\cal O}_m={\cal O}((m,0);2\nu)$.
For a fixed flat $F\in{\cal O}_i$, let ${\cal O}_j'(F)$ be the set
of all flats in  ${\cal O}_j$ containing $F$, where $0\leq i\leq j\leq \nu$.

 \begin{lemma}\label{lemma2.4}{\rm(See \cite[Lemma~9.4.1]{Brouwer} and \cite[Theorems~1.17,~1.19,~3.38,~5.37,~6.43]{Wan}).}
 Let $0\leq i\leq j\leq \nu$. Then the following hold.
 \begin{itemize}
 \item[\rm(i)]
 The size of ${\cal O}_i$ is $q^{2\nu-i}{\nu\brack i}_q\prod_{t=\nu-i+1}^{\nu}(q^{t+e-1}+1)$.

 \item[\rm(ii)]
 For a fixed flat $F\in{\cal O}_i$, the size of ${\cal O}_j'(F)$ is ${\nu-i\brack j-i}_q\prod_{t=\nu-j+1}^{\nu-i}(q^{t+e-1}+1)$.
 \end{itemize}
 \end{lemma}

Recall that $M$ is the incidence matrix with rows indexed with vectors in $\mathbb{F}_q^{2\nu}$
and columns indexed with flats in ${\cal O}_{\nu}$ such that entry $M_{x,F}=1$ if
and only if $x\in F$.

\begin{lemma}\label{lemma2.5}
If $\nu\geq1$, then the rank of $M$ is $(q^{\nu}-1)(q^{\nu-1+e}+1)+1$ over the real field $\mathbb{R}$.
\end{lemma}
\proof Let $N=MM^{\top}$. Then both $N$ and
$M$ have the same rank over $\mathbb{R}$. Note that $N$ is a $q^{2\nu}\times q^{2\nu}$ matrix with rows and columns indexed with the vectors in $\mathbb{F}_q^{2\nu}$.
 For any $x,y\in\mathbb{F}_q^{2\nu}$, the entry
$N_{xy}$ is the number of maximal totally isotropic flats containing both $x$ and $y$.
By Lemma~\ref{lemma2.4}, we have
$$N_{xy}=\left\{\begin{array}{ll}
\prod_{t=1}^{\nu}(q^{t+e-1}+1)  &  \hbox{if}\;x=y,\\
\prod_{t=1}^{\nu-1}(q^{t+e-1}+1)   & \hbox{if}\;\{x\}\vee\{y\}\in{\cal O}_1,\\
0           &   \hbox{otherwise},
\end{array}\right.$$
which implies that
\begin{equation}\label{equa2.6}
N=I\prod_{t=1}^{\nu}(q^{t+e-1}+1)+A\prod_{t=1}^{\nu-1}(q^{t+e-1}+1),
\end{equation}
where $A$ is the adjacency matrix of the graph ${\cal G}$.

For the symplectic case,  by Lemma~\ref{lemma2.3} and (\ref{equa2.6}), $N$ is nonsingular, which implies that the rank of  $N$ is  $q^{2\nu}$.
For the unitary and orthogonal cases,  by Lemma~\ref{lemma2.3} and (\ref{equa2.6}) again,
the distinct eigenvalues of $N$ are $q^{\nu}\prod_{t=1}^{\nu}(q^{t+e-1}+1),q^{\nu}\prod_{t=1}^{\nu-1}(q^{t+e-1}+1)$ and $0$, and the corresponding multiplicities are $1,(q^{\nu}-1)(q^{\nu+e-1}+1)$ and
$(q^{\nu}-1)(q^{\nu}-q^{\nu+e-1})$, respectively.
So, the rank of $N$ is $(q^{\nu}-1)(q^{\nu-1+e}+1)+1$. $\qed$

\section{Association schemes}
Throughout this paper, we will denote the all-one  matrix by $J$. The all-one column vector of length
$n$ will be denoted by ${\rm j}_n$. For convenience, we write ${\rm j}_n$ by ${\rm j}$ if its length is clear from the context.

A $d$-{\it class association scheme} ${\mathfrak X}$ is a pair $(X,\{R_i\}_{i=0}^{d})$, where $X$ is a finite set, and each $R_{i}$ is a nonempty subset of $X\times X$ satisfying the following axioms:
\begin{itemize}

\item[\rm(i)]
$R_0=\{(x,x) : x\in X\}$;

\item[\rm(ii)]
$\{R_i\}_{0\leq i\leq d}$ is a partition of $X\times X$;

\item[\rm(iii)]
${}^tR_i=R_{i}$ for each $i\in\{0,1,\ldots,d\}$, where
${}^tR_i=\{(y,x) : (x,y)\in R_i\}$;

\item [\rm(iv)]
there exist integers $p_{ij}^{k}$ such that for all $(x,y)\in R_{k}$ and all $i,j,k$,
$$
p_{ij}^{k}=|\{z\in X : (x,z)\in
R_i,(z,y)\in R_j\}|.
$$
\end{itemize}
The integers $p_{ij}^{k}$ are called the {\em intersection numbers}
of ${\mathfrak X}$, and $v_{i}\,(=p_{ii}^{0})$ is called the
{\it valency} of $R_i$.

Let $M_X(\mathbb{R})$ be the algebra of matrices over the real field $\mathbb{R}$ with rows and columns indexed by $X$. The
$i$th adjacency matrix $A_i\in M_X(\mathbb{R})$ of $\mathfrak{X}$ is defined to be the matrix of degree $|X|$ in which
$(x,y)$-entries are equal to 1 if and only if $(x,y)\in R_i$. Let
$\mathfrak{A} = \langle A_0,A_1,\ldots,A_d\rangle$ be the subalgebra of $M_X(\mathbb{R})$
spanned by $A_0,A_1,\ldots,A_d$. Then $\mathfrak{A}$ is
called the {\it Bose-Mesner algebra} of $\mathfrak{X}$.

Let $V_0\perp V_1\perp\cdots\perp V_d$ be the orthogonal decomposition  of $\mathbb{R}^{X}$ in common
eigenspaces of $A_0,A_1,\ldots,A_d$. It
is well known that $\mathfrak{A}$ is semi-simple, hence $\mathfrak{A}$ has a basis consisting of the primitive
idempotents $E_0,E_1,\ldots,E_d$. The {\it first eigenmatrix} $P=(p_i(j))$ and
the {\it second eigenmatrix} $Q=(q_i(j))$ of $\mathfrak X$ are defined by
$A_i=\sum_{j=0}^dp_i(j)E_j$ and
$E_i=\frac{1}{|X|}\sum_{j=0}^dq_i(j)A_j$, respectively.
Note that $P$ and $Q$ satisfy $PQ=QP=|X|I$, and it is well known that
\begin{equation}\label{equa3.1}
m_j=\frac{|X|}{\sum_{t=0}^d\frac{(p_t(j))^2}{v_t}}\quad\hbox{and}\quad \frac{q_j(i)}{m_j}=\frac{p_i(j)}{v_i}
\end{equation}
where $m_j={\rm Tr}(E_j)$ is called the $j$th
{\it multiplicity}. The theory of association schemes was discussed in \cite{Bannai, Brouwer}.

From $J=\sum_{i=0}^{d}A_i$, we deduce that $\sum_{i=0}^dp_i(j),\,0\leq j\leq d$,
are eigenvalues of $J$. Note that the distinct eigenvalues of $J$ are $|X|$ and $0$,
and the corresponding multiplicities are $1$ and $|X|-1$, respectively. Therefore, we obtain
\begin{equation}\label{equa3.1.1}
\sum_{i=0}^{d}p_{i}(j)=\left\{\begin{array}{ll}
|X| & \hbox{if}\;j=0,\\
0   &\hbox{otherwise}.
\end{array}\right.
\end{equation}

Let ${\mathfrak X}=(X,\{R_i\}_{i=0}^{d})$ be a $d$-class
association scheme. Let $S$ be a subset of $X$ and $\chi_S$ be its characteristic
vector, i.e. $\chi_S$ is the column vector in which its positions correspond to the elements of $X$, such that $(\chi_{{\cal S}})_{x}=1$ if $x\in S$ and 0 otherwise.
The row vector $u=(u_0,u_1,\ldots,u_d)\in\mathbb{R}^{d+1}$ with
$$u_i=\frac{|R_i\cap(S\times S)|}{|S|}=\frac{\chi_S^{\top}A_i\chi_S}{|S|}$$
is called the {\it inner distribution} of $S$.

\begin{lemma}\label{lemma3.1}{\rm (See \cite[Lemma~2.5.1 and Proposition~2.5.2]{Brouwer}).}
Let ${\mathfrak X}=(X,\{R_i\}_{i=0}^{d})$ be a $d$-class
association scheme, and $\mathfrak{A}$ be its Bose-Mesner algebra.
 Denote the primitive idempotent basis
of $\mathfrak{A}$ by $\{E_i : 0\leq i\leq d\}$, with common eigenspaces $V_0,\ldots,V_d$. Let $S$ be a subset of $X$ and $\chi_S$ be its characteristic
vector. Then there exist $v_i\in V_i$ and $a_i\in\mathbb{R}$  for $i=0,1,\ldots,d$, such that
$\chi_S=\sum_{i=0}^da_iv_i$. Moreover, if $u$ is the inner distribution of $S$, then
the following properties are equivalent for fixed $0\leq i\leq d$.
\begin{itemize}
\item[\rm(i)]
$(uQ)_i=0$, where $Q$ is the  second eigenmatrix  of $\mathfrak X$.

\item[\rm(ii)]
$E_i\chi_S=0$.
\end{itemize}
The $E_i\chi_S=0$ implies that the projection of $\chi_S$ onto the eigenspace $V_i$ is zero, thus
$a_i=0$.
\end{lemma}

\begin{lemma}\label{lemma3.1.2}{\rm (See \cite[Theorem~6.8]{Delsarte}).}
Let ${\mathfrak X}=(X,\{R_i\}_{i=0}^{d})$ be a $d$-class
association scheme, with  $\{E_i : 0\leq i\leq d\}$ the primitive idempotent basis
 of the Bose-Mesner algebra of ${\mathfrak X}$. Suppose $G$ is a subgroup
of {\rm Aut($\mathfrak{X}$)} that acts transitively on $X$ and its orbits on $X\times X$ are the relations $R_0,R_1, \ldots, R_d$.
Let $\chi$ and $\psi$ be vectors of $\mathbb{R}^{|X|}$. Then the following two statements are equivalent.
\begin{itemize}
\item[\rm(i)]
For all $i\geq1$, we have $E_i\chi=0$ or $E_i\psi=0$.

\item[\rm(ii)]
$\chi\cdot\psi^g$ is constant for all $g\in G$.
\end{itemize}
\end{lemma}

Let $\mathbb{F}_q^{2\nu}$ be the $2\nu$-dimensional classical space with parameter $e$, and ${\cal M}_{\nu}$ be the set of all maximal totally isotropic subspaces in $\mathbb{F}_q^{2\nu}$. For $0\leq i\leq \nu$, define
$$R_i=\{(P,Q)\in{\cal M}_{\nu}\times{\cal M}_{\nu} : \dim(P\cap Q)=\nu-i\}.$$
Then ${\mathfrak
X}({\cal M}_{\nu}):=({\cal M}_{\nu},\{R_{i}\}_{0\leq i\leq \nu})$ is the {\it dual polar scheme} of
$\mathbb{F}_q^{2\nu}$. Let $A_i$ be the $i$th adjacency matrix of ${\mathfrak X}({\cal M}_{\nu})$. Then there is an orthogonal decomposition $V_0\perp V_1\perp\cdots\perp V_\nu$ of $\mathbb{R}^{{\cal M}_{\nu}}$ in common
eigenspaces of $A_0,A_1,\ldots,A_\nu$.

\begin{lemma}\label{lemma3.2}{\rm (See \cite[Theorem~9.4.3]{Brouwer} and \cite[Theorem~4.3.6]{Vanhove}).}
Let ${\mathfrak X}({\cal M}_{\nu})$ be the dual polar scheme of
$\mathbb{F}_q^{2\nu}$, and $(p^{(2\nu)}_{i}(j))$ be the first eigenmatrix of ${\mathfrak X}({\cal M}_{\nu})$. Then the eigenvalue of $A_i$ corresponding to $V_j$ is
$$p^{(2\nu)}_i(j)=\sum_{s=\max\{0,j-i\}}^{\min\{j,\nu-i\}}(-1)^{j+s}{j\brack s}_q{\nu-j\brack \nu-i-s}_qq^{e(i+s-j)+{j-s\choose 2}+{i+s-j\choose 2}},$$
and $V_j$ has dimension
$$m_j^{(2\nu)}=q^{j}{\nu\brack j}_q\frac{q^{\nu+e-2j}+1}{q^{\nu+e-j}+1}\prod_{s=1}^j\frac{q^{\nu+e-s}+1}{q^{s-e}+1}.$$
In particular,
$$v^{(2\nu)}_i=p^{(2\nu)}_i(0)=q^{i(i+2e-1)/2}{\nu\brack i}_q\quad\hbox{and}\quad p^{(2\nu)}_1(j)=q^e{\nu-j\brack 1}_q-{j\brack 1}_q.$$
\end{lemma}

For $0\leq i\leq\nu$, define
\begin{eqnarray*}
R_{(i,0)}&=&\{(P+x,Q+y)\in {\cal O}_{\nu}\times{\cal O}_{\nu} : (P+x)\cap (Q+y)\not=\emptyset, \dim(P\cap Q)=\nu-i\},\\
R_{(i,1)}&=&\{(P+x,Q+y)\in {\cal O}_{\nu}\times{\cal O}_{\nu} : (P+x)\cap (Q+y)=\emptyset, \dim(P\cap Q)=\nu-i\}.
\end{eqnarray*}
By \cite[Theorem~1.1]{Guo},  the configuration
\begin{equation}\label{scheme}
{\mathfrak X}({\cal O}_{\nu}):=({\cal O}_{\nu},\{R_{(0,0)},R_{(0,1)},\dots,R_{(\nu-1,0)},R_{(\nu-1,1)},R_{(\nu,0)}\})
\end{equation} is an association scheme   and called the {\it maximal totally isotropic flat  scheme}
in $ACG(2\nu,\mathbb{F}_q)$. For the symplectic and orthogonal cases, the scheme ${\mathfrak X}({\cal O}_{\nu})$
was also obtained in \cite{Liu,Liu2}.
Let $A_{(i,\xi)}$ be the $(i,\xi)$th adjacency matrix of ${\mathfrak X}({\cal O}_{\nu})$. Then there is an orthogonal decomposition $V_{(0,0)}\perp V_{(0,1)}\perp\cdots\perp V_{(\nu-1,0)}\perp V_{(\nu-1,1)}\perp V_{(\nu,0)}$ of $\mathbb{R}^{{\cal O}_{\nu}}$ in common
eigenspaces of $A_{(0,0)},A_{(0,1)},\ldots,A_{(\nu-1,0)},A_{(\nu-1,1)},A_{(\nu,0)}$.

\begin{lemma}\label{lemma-trs}{\rm (See \cite[Lemma~3.1]{Guo}).}
The group $AG_{2\nu}$ is a subgroup of Aut(${\mathfrak X}({\cal O}_{\nu})$) that acts transitively on ${\cal O}_{\nu}$ and its orbits on ${\cal O}_{\nu}\times{\cal O}_{\nu}$
are exactly the relations $$R_{(0,0)},R_{(0,1)},\dots,R_{(\nu-1,0)},R_{(\nu-1,1)},R_{(\nu,0)}.$$
\end{lemma}

\begin{lemma}\label{lemma3.3}
Let ${\mathfrak X}({\cal O}_{\nu})$ be as in (\ref{scheme}), and $v_{(i,\xi)}$ be the valency of the relation $R_{(i,\xi)}$. Then
$$v_{(i,0)}=q^{i(i+2e+1)/2}{\nu\brack i}_q\quad\hbox{and}\quad
v_{(i,1)}=(q^{\nu-i}-1)q^{i(i+2e+1)/2}{\nu\brack i}_q.$$
\end{lemma}
\proof By Lemma~\ref{lemma3.2}, we obtain that  the valency of the dual polar scheme ${\mathfrak X}({\cal M}_{\nu})$ of
$\mathbb{F}_q^{2\nu}$ is $v^{(2\nu)}_i$. By Lemma~\ref{lemma2.1},
$v_{(i,0)}=q^iv^{(2\nu)}_i$ and $v_{(i,1)}=(q^{\nu}-q^i)v^{(2\nu)}_i$. Therefore, the desired result follows. $\qed$

\begin{lemma}\label{lemma3.4}{\rm (See \cite[Theorem~1.1]{Guo}).}
Let ${\mathfrak X}({\cal O}_{\nu})$ be as in (\ref{scheme}), and $(p_{(i,\xi)}(j,\eta))$ be the first eigenmatrix of ${\mathfrak X}({\cal O}_{\nu})$.
Suppose that $(p^{(2\nu)}_i(j))$ be the first eigenmatrix of ${\mathfrak X}({\cal M}_{\nu})$, and let
 $$C^{(\ell)}=(c^{(\ell)}_{\xi}(\eta))
 =\left(\begin{array}{cc}
 1 & q^{\ell}-1\\
 1  & -1
 \end{array}\right)$$
 be the first eigenmatrix of the complete graph of size $q^{\ell}$.
Then the eigenvalue of $A_{(i,\xi)}$ corresponding to $V_{(j,\eta)}$ is
$$p_{(i,\xi)}(j,\eta)=q^ic^{(\nu-i)}_{\xi}(\eta)p^{(2\nu-2\eta)}_i(j),$$
and $V_{(j,\eta)}$ has dimension
$$
m_{(j,\eta)}=\left\{\begin{array}{ll}
q^{j}{\nu\brack j}_q\frac{q^{\nu+e-2j}+1}{q^{\nu+e-j}+1}\prod_{s=1}^j\frac{q^{\nu+e-s}+1}{q^{s-e}+1}, &\hbox{if}\;\eta=0,\\
(q^{\nu}-1)(q^{\nu+e-1}+1)q^{j}{\nu-1\brack j}_q\frac{q^{\nu-1+e-2j}+1}{q^{\nu-1+e-j}+1}\prod_{s=1}^j\frac{q^{\nu-1+e-s}+1}{q^{s-e}+1} &\hbox{if}\;\eta=1.
\end{array}\right.$$
\end{lemma}
\proof By \cite[Theorem~1.1]{Guo}, the eigenvalue $p_{(i,\xi)}(j,\eta)$ of $A_{(i,\xi)}$ corresponding to $V_{(j,\eta)}$ is computed. Let $m_{(j,\eta)}$ be the dimension of $V_{(j,\eta)}$. By (\ref{equa3.1}) and Lemma~\ref{lemma3.2}, we obtain
\begin{eqnarray}
m_j^{(2\nu-2\eta)}&=&\frac{|{\cal M}_{\nu-\eta}|}{\sum\limits_{t=0}^{\nu-\eta}\frac{\left(p^{(2\nu-2\eta)}_{t}(j)\right)^2}{v^{(2\nu-2\eta)}_{t}}}
=q^{j}{\nu-\eta\brack j}_q\frac{q^{\nu-\eta+e-2j}+1}{q^{\nu-\eta+e-j}+1}\prod_{s=1}^j\frac{q^{\nu-\eta+e-s}+1}{q^{s-e}+1},\label{equa3.2}\\
m_{(j,\eta)}&=&\frac{|{\cal O}_{\nu}|}{\sum\limits_{(t,\zeta)}\frac{\left(p_{(t,\zeta)}(j,\eta)\right)^2}{v_{(t,\zeta)}}}
=\frac{q^{\nu}|{\cal M}_{\nu}|}{\sum\limits_{(t,\zeta)}\frac{\left(q^tc^{(\nu-t)}_{\zeta}(\eta)p^{(2\nu-2\eta)}_t(j)\right)^2}{v_{(t,\zeta)}}}.\label{equa3.3}
\end{eqnarray}
So, by Lemmas~\ref{lemma3.2}, \ref{lemma3.3}, (\ref{equa3.2}) and (\ref{equa3.3}), we have
\begin{eqnarray*}
m_{(j,0)}&=&\frac{q^{\nu}|{\cal M}_{\nu}|}{\sum\limits_{t=0}^{\nu}\frac{q^t\left(p^{(2\nu)}_t(j)\right)^2}{v^{(2\nu)}_{(t)}}
+\sum\limits_{t=0}^{\nu}\frac{(q^{\nu}-q^t)\left(p^{(2\nu)}_t(j)\right)^2}{v^{(2\nu)}_{(t)}}}
=\frac{q^{\nu}|{\cal M}_{\nu}|}{q^{\nu}\sum\limits_{t=0}^{\nu}\frac{\left(p^{(2\nu)}_t(j)\right)^2}{v^{(2\nu)}_{(t)}}}
=m_j^{(2\nu)},\\
m_{(j,1)}&=&\frac{q^{\nu}|{\cal M}_{\nu}|}{\sum\limits_{t=0}^{\nu}\frac{q^t\left(p^{(2\nu-2)}_t(j)\right)^2}{v^{(2\nu)}_{(t)}}
+\sum\limits_{t=0}^{\nu-1}\frac{q^t\left(p^{(2\nu-2)}_t(j)\right)^2}{(q^{\nu-t}-1)v^{(2\nu)}_{(t)}}}\\
&=&\frac{q^{\nu}(q^{\nu+e-1}+1)|{\cal M}_{\nu-1}|}{{\sum\limits_{t=0}^{\nu-1}\frac{(q^{\nu}-q^t)\left(p^{(2\nu-2)}_t(j)\right)^2}{(q^{\nu}-1)v^{(2\nu-2)}_{(t)}}
+\sum\limits_{t=0}^{\nu-1}\frac{q^t\left(p^{(2\nu-2)}_t(j)\right)^2}{(q^{\nu}-1)v^{(2\nu-2)}_{(t)}}}}\\
&=&\frac{(q^{\nu}-1)(q^{\nu+e-1}+1)|{\cal M}_{\nu-1}|}{{\sum\limits_{t=0}^{\nu-1}\frac{\left(p^{(2\nu-2)}_t(j)\right)^2}{v^{(2\nu-2)}_{(t)}}}}
=(q^{\nu}-1)(q^{\nu+e-1}+1)m_j^{(2\nu-2)}.
\end{eqnarray*}
Therefore, the desired result follows. $\qed$

Now we introduce $\phi^{(2\nu)}_{i}(j)=\max\{k : q^k | p^{(2\nu)}_{i}(j)\}$,
the exponent of $q$ in $p^{(2\nu)}_{i}(j)$. If $p^{(2\nu)}_{i}(j)=0$, we put $\phi^{(2\nu)}_{i}(j)=\infty$.
Then \begin{equation}\label{equ-01}
\phi^{(2\nu)}_{i}(0)={i\choose 2}+ei \quad\hbox{and}\quad \phi^{(2\nu)}_{i}(1)={i-1\choose 2}+e(i-1).
\end{equation} For $j\geq 2$, we obtain the following result.

\begin{lemma}\label{lemma3.7-N}
Let $2\leq i,j\leq\nu$. Then
$\phi^{(2\nu)}_{i}(j)={i\choose 2}+(j-i)(j-e)$ if $j-\frac{i}{2}-\frac{e}{2}<0$,
 $\phi^{(2\nu)}_{i}(j)$ is as in Table~1 if $0\leq j-\frac{i}{2}-\frac{e}{2}\leq\nu-i$,
and $\phi^{(2\nu)}_{i}(j)=
(j-e-\nu+1)(j-\nu+i-1)+{i-1\choose 2}+e(i-1)$ if $j-\frac{i}{2}-\frac{e}{2}>\nu-i$.
\end{lemma}
\proof By  \cite[Lemma~2.4]{De Boeck0}, $\phi^{(2\nu)}_{i}(j)$ is computed for all $2\leq i,j\leq\nu$. $\qed$

\begin{table} \caption{$0\leq j-\frac{i}{2}-\frac{e}{2}\leq\nu-i$.}
 \begin{center}
\begin{tabular}{lll}
\hline Type & $i$ &  $\phi^{(2\nu)}_{i}(j)$  \\
\hline Orthogonal & even & $\frac{i(i-2)}{4}$   \\
Orthogonal & odd & $\left\{\begin{array}{ll}
                       \frac{(i-1)^2}{4} & \hbox{if}\;j\not=\frac{\nu}{2}\\
                       \infty  & \hbox{if}\;j=\frac{\nu}{2}
                       \end{array}\right.$   \\
 Unitary & all & $\frac{i(i-1)}{4}$    \\
  Symplectic & even  & $\left\{\begin{array}{ll}
                         \infty &\hbox{if}\;2(j-1)=2i=\nu\equiv 0 \;(4)\\
                         \frac{i^2}{4} & \hbox{otherwise}
                         \end{array}\right.$  \\
  Symplectic & odd & $\frac{i^2-1}{4}$  \\
   \hline
\end{tabular}
 \end{center}
\end{table}

\begin{lemma}\label{lemma3.8-N}
Let $2\leq i\leq\nu$. Then $\phi^{(2\nu)}_{i}(0)\not=\phi^{(2\nu)}_{i}(j)$ for all $1\leq j\leq \nu$, except $(j,e)=(\nu,0)$.
Moreover, if $e=0$, then $p^{(2\nu)}_{i}(0)=p^{(2\nu)}_{i}(\nu)$ if and only if $i$ is even.
\end{lemma}
\proof
By (\ref{equ-01}), $\phi^{(2\nu)}_{i}(0)\not=\phi^{(2\nu)}_{i}(1)$.
If $2\leq j<\frac{i}{2}+\frac{e}{2}$, by Lemma~\ref{lemma3.7-N} and (\ref{equ-01}),
$$\phi^{(2\nu)}_{i}(0)=\phi^{(2\nu)}_{i}(j)\Leftrightarrow ei=(j-i)(j-e)\Leftrightarrow j-i=e,$$
a contradiction since $j<i+e$.
If $\max\{2,\frac{i}{2}+\frac{e}{2}\}\leq j\leq \nu-\frac{i}{2}+\frac{e}{2}$, by Table~1 in Lemma~\ref{lemma3.7-N} and (\ref{equ-01}),
it is easy to prove $\phi^{(2\nu)}_{i}(0)\not=\phi^{(2\nu)}_{i}(j)$.
If $j\geq\max\{2,\nu-\frac{i}{2}+\frac{e}{2}+1\}$, by Lemma~\ref{lemma3.7-N} and (\ref{equ-01}),
\begin{eqnarray*}
\phi^{(2\nu)}_{i}(0)=\phi^{(2\nu)}_{i}(j)&\Leftrightarrow& i-1+e=(j-e-\nu+1)(j-\nu+i-1)\\
& \Leftrightarrow& (j-\nu+i)(j-\nu-e)=0 \\
&\Leftrightarrow& j=\nu-i\;\hbox{or}\;j=\nu+e.
\end{eqnarray*}
If $j=\nu-i$, a contradiction since $j>\nu-\frac{i}{2}+\frac{e}{2}$. If $j=\nu+e$ then $(j,e)=(\nu,0)$, which implies that
 $\phi^{(2\nu)}_{i}(0)=\phi^{(2\nu)}_{i}(j)$ if and only if $(j,e)=(\nu,0)$.
 If $e=0$, by Lemma~\ref{lemma3.2} and (\ref{equ-01}),
 $p^{(2\nu)}_i(\nu)=(-1)^{i}q^{{i\choose 2}}{\nu\brack i}_q=p^{(2\nu)}_i(0)$ if and only if $i$ is even. $\qed$

\begin{lemma}\label{lemma3.9-N}
Let ${\mathfrak X}({\cal O}_{\nu})$ be as in (\ref{scheme}), and  $(p_{(i,\xi)}(j,\eta))$ be the first eigenmatrix of ${\mathfrak X}({\cal O}_{\nu})$. Suppose $(i,\xi)\not=(0,0)$. Then the eigenvalue $p_{(i,\xi)}(0,1)$ of $A_{(i,\xi)}$ only corresponds to $V_{(0,1)}$, except in the following cases: {\rm(a)} $\nu\geq2$ and $(i,\xi)=(0,1)$; {\rm(b)} $\nu\geq2$ and $(i,\xi)=(\nu,0)$; {\rm(c)}  $i$ is even and $2\leq i\leq\nu-1$.
\end{lemma}
\proof We divide into the following five cases.

{\it Case}~1: $(i,\xi)=(0,1)$. If $\nu=1$,
by $p_{(0,1)}(0,1)=-1\not=1=p_{(0,1)}(0,0)$ and $p_{(0,1)}(0,1)=-1\not=q-1=p_{(0,1)}(1,0)$,
the eigenvalue $p_{(0,1)}(0,1)$ of $A_{(0,1)}$  only corresponds to $V_{(0,1)}$. If $\nu\geq2$,
 by Lemmas~\ref{lemma3.2} and~\ref{lemma3.4},
$$p_{(0,1)}(j,\eta)=\left\{\begin{array}{ll}
 q^{\nu}-1     & \hbox{if}\;\eta=0,\\
 -1     & \hbox{if}\;\eta=1,
\end{array}\right.$$
which implies that $p_{(0,1)}(0,1)=p_{(0,1)}(j,1)$ for all $1\leq j\leq \nu-1$.

{\it Case}~2: $(i,\xi)=(1,0)$. If $\nu=1$, then
$$p_{(1,0)}(0,1)=0\not=1=p_{(1,0)}(0,0)\quad \hbox{and}\quad p_{(1,0)}(0,1)=0\not=-q=p_{(1,0)}(1,0).$$
If $\nu\geq2$,  by Lemmas~\ref{lemma3.2} and~\ref{lemma3.4},
 $$p_{(1,0)}(0,1)=q^{1+e}{\nu-1\brack 1}_q\not=q\left(q^{e}{\nu-1-j\brack 1}_q-{j\brack 1}_q\right)=p_{(1,0)}(j,1)$$
  for all $1\leq j\leq \nu-1$, and $$p_{(1,0)}(0,1)=q^{1+e}{\nu-1\brack 1}_q\not=q\left(q^e{\nu-j\brack 1}_q-{j\brack 1}_q\right)=p_{(1,0)}(j,0)$$
   for all $0\leq j\leq\nu$.
 So, the eigenvalue $p_{(0,1)}(0,1)$ of $A_{(0,1)}$  only corresponds to $V_{(0,1)}$.

{\it Case}~3: $(i,\xi)=(1,1)$. Then $\nu\geq2$. By Lemmas~\ref{lemma3.2} and~\ref{lemma3.4},
$$p_{(1,1)}(0,1)=-q^{1+e}{\nu-1\brack 1}_q\not=-q\left(q^{e}{\nu-1-j\brack 1}_q-{j\brack 1}_q\right)=p_{(1,1)}(j,1)$$ for all $1\leq j\leq \nu-1$. Suppose $$p_{(1,1)}(0,1)=-q^{1+e}{\nu-1\brack 1}_q=q(q^{\nu-1}-1)\left(q^e{\nu-j\brack 1}_q-{j\brack 1}_q\right)=p_{(1,1)}(j,0).$$
Then $q{\nu-1\brack 1}_q(q^{e+\nu-j}-q^j+1)=0$, which implies that $q^{e+\nu-j}-q^j+1=0$.
If $j=0$, then $q^{\nu+e}=0$, a contradiction. If  $j=\nu$, then $q^{\nu}-q^{e}=1$, a contradiction since $q^{\nu}-q^{e}\geq q^{2}-q^{e}>1$.
If $1\leq j<\nu$, then $q^j-q^{\nu-j+e}=1$,  a contradiction since $q^j-q^{\nu-j+e}$ is divisible by $q$, but $1$ is not divisible by $q$.
So, the eigenvalue $p_{(1,1)}(0,1)$ of $A_{(1,1)}$  only corresponds to $V_{(0,1)}$.

{\it Case}~4: $(i,\xi)=(\nu,0)$ and $\nu\geq2$. Then $p_{(\nu,0)}(0,1)=p_{(\nu,0)}(j,1)=0$ for all $1\leq j\leq \nu-1$. So,
there exists a $j=1$ such that $p_{(\nu,0)}(0,1)=p_{(\nu,0)}(1,1)$.

{\it Case}~5: $2\leq i\leq\nu-1$.
Let $\phi_{(i,\xi)}(j,\eta)=\max\{k : q^k | p_{(i,\xi)}(j,\eta)\}$,
the exponent of $q$ in $p_{(i,\xi)}(j,\eta)$. If $p_{(i,\xi)}(j,\eta)=0$, we  put $\phi_{(i,\xi)}(j,\eta)=\infty$.
Note that $\phi_{(i,\xi)}(0,1)\not=\phi_{(i,\xi)}(j,\eta)$ implies $p_{(i,\xi)}(0,1)\not=p_{(i,\xi)}(j,\eta)$.
 By Lemmas~\ref{lemma3.2} and~\ref{lemma3.4}, we have
\begin{equation}\label{equ-i}
\phi_{(i,\xi)}(j,\eta)=\left\{\begin{array}{ll}
i+\phi^{(2\nu)}_i(j) & \hbox{if}\;(\xi,\eta)=(0,0),\\
i+\phi^{(2\nu-2)}_i(j)    & \hbox{if}\;(\xi,\eta)=(0,1),\\
i+\phi^{(2\nu)}_i(j)      & \hbox{if}\;(\xi,\eta)=(1,0),\\
i+ \phi^{(2\nu-2)}_i(j)      & \hbox{if}\;(\xi,\eta)=(1,1).
\end{array}\right.
\end{equation}
By Lemma~\ref{lemma3.8-N},
$\phi^{(2\nu-2)}_i(0)\not=\phi^{(2\nu-2)}_i(j)$ for all $1\leq j\leq\nu-1$, except $(j,e)=(\nu-1,0)$. So, from (\ref{equ-i}) we deduce that
 $\phi_{(i,\xi)}(0,1)\not=\phi_{(i,\xi)}(j,1)$ for all $1\leq j\leq\nu-1$, except $(j,e)=(\nu-1,0)$.
 If $(j,e)=(\nu-1,0)$, then
 $$p_{(i,\xi)}(0,1)=p_{(i,\xi)}(\nu,1)\Leftrightarrow
 q^{i+{i\choose 2}}{\nu-1\brack i}_q=(-1)^{i}q^{i+{i\choose 2}}{\nu-1\brack i}_q \Leftrightarrow i\;\hbox{even}
 $$

 Next, we compare $p_{(i,\xi)}(0,1)$ and $p_{(i,\xi)}(j,0)$ for all $0\leq j\leq\nu$.
 It is easy to show $p_{(i,\xi)}(0,1)\not=p_{(i,\xi)}(0,0)$.
 Since $i\leq\nu-1$, by Lemma~\ref{lemma3.7-N} and (\ref{equ-i}), $$\phi_{(i,\xi)}(0,1)=i+\phi^{(2\nu-2)}_{i}(0)=i+\phi^{(2\nu)}_{i}(0).$$
  By Lemma~\ref{lemma3.8-N}, $\phi_{(i,\xi)}(0,1)=i+\phi^{(2\nu)}_{i}(0)\not=i+\phi^{(2\nu)}_{i}(j)=\phi_{(i,\xi)}(j,0)$ for all $1\leq j\leq \nu$, except $(j,e)=(\nu,0)$. If $e=0$, since ${\nu\brack i}_q>{\nu-1\brack i}_q$,
  $$p_{(i,\xi)}(0,1)=(-1)^{\xi}q^{i+{i\choose 2}}{\nu-1\brack i}_q
  \not=p_{(i,\xi)}(\nu,0)  =\left\{\begin{array}{ll}
  (-1)^{i}q^{i+{i\choose 2}}{\nu\brack i}_q  & \hbox{if}\;\xi=0,\\
  (-1)^{i}q^{i+{i\choose 2}}(q^{\nu-i}-1){\nu\brack i}_q  & \hbox{if}\;\xi=1.
  \end{array}\right.$$

 From above the discussion, we obtain that the desired result follows. $\qed$

\section{Equivalent definitions}
In this section,   we give several equivalent definitions for a Cameron-Liebler set in ${\cal O}_{\nu}$, see Definition~\ref{defin1.1}.
In the rest of this paper, we always assume that ${\mathfrak X}({\cal O}_{\nu})$ is as in (\ref{scheme}), $\mathfrak{A}=\langle A_{(0,0)},A_{(0,1)},\dots,$ $A_{(\nu-1,0)},A_{(\nu-1,1)},A_{(\nu,0)}\rangle$ is its Bose-Mesner algebra,
 and $(p_{(i,\xi)}(j,\eta))$ and $(q_{(i,\xi)}(j,\eta))$ are the first eigenmatrix and the second eigenmatrix of ${\mathfrak X}({\cal O}_{\nu})$, respectively.
 Denote the primitive idempotent basis
of $\mathfrak{A}$ by $\{E_{(0,0)},E_{(0,1)},\dots,E_{(\nu-1,0)},E_{(\nu-1,1)},E_{(\nu,0)}\}$, with common eigenspaces
$V_{(0,0)},V_{(0,1)},\dots, V_{(\nu-1,0)},V_{(\nu-1,1)},V_{(\nu,0)}.$

Recall that $ACG(2\nu,\mathbb{F}_q)$ is the $2\nu$-dimensional classical affine
space with parameter $e$, and ${\cal O}_\nu$ is the set of all maximal totally isotropic flats in $ACG(2\nu,\mathbb{F}_q)$.

\begin{defin} Let $1\leq i\leq\nu$, and $T$ be a  $(\nu+i,2i)$-flat  in $ACG(2\nu,\,\mathbb{F}_q)$.
\begin{itemize}
\item[\rm(i)]
A partial $(\nu,0)$-spread in $T$ is a set of pairwise disjoint $(\nu,0)$-flats contained in $T$.

\item[\rm(ii)]
A $(\nu,0)$-conjugated switching set in $T$ is a pair of disjoint partial $(\nu,0)$-spreads in $T$ that cover the same
set of vectors in $T$.

\item[\rm(iii)]
A $(\nu,0)$-spread in $T$ is a partial $(\nu,0)$-spread in $T$ that partitions the set of all vectors in $T$.
\end{itemize}
\end{defin}

Note that $\mathbb{F}_q^{2\nu}$ is a $(2\nu,2\nu)$-flat in $ACG(2\nu,\mathbb{F}_q)$.
Next, we give some examples of $(\nu,0)$-spreads in $\mathbb{F}_q^{2\nu}$.

\begin{lemma}\label{lemma4.1}
Let $\nu\geq 1$. Then the following sets ${\cal S}$
are $(\nu,0)$-spreads in $\mathbb{F}_q^{2\nu}$.
\begin{itemize}
\item[\rm(i)] {\rm (Type I)} Define ${\cal S}$ as the set of all the cosets of
$\mathbb{F}_q^{2\nu}$ relative to the maximal totally isotropic subspace $P$. That is, ${\cal S}=\{P+x: x\in\mathbb{F}_q^{2\nu}\}$.

\item[\rm(ii)] {\rm (Type I\!I)}
Suppose that $\nu\geq2$ and $Q$ is a subspace of type $(\nu+1,2)$
in $\mathbb{F}_q^{2\nu}$. Let $P_1$ and $P_2$ be two distinct maximal totally isotropic subspaces contained in $Q$.
Define $${\cal S}=\{P_1+y_1 : y_1\in Q \}\cup\{P_2+y_2 : y_2\in \mathbb{F}_q^{2\nu}\setminus Q\}.$$
\end{itemize}
\end{lemma}
\proof (i) Trivial.

(ii) By \cite[Theorems~3.27, 5.28 and 6.33]{Wan}, the number of maximal totally isotropic subspaces contained in $Q$ is
$q^e+1\geq2$, which implies that $Q$ contains two distinct maximal totally isotropic subspaces $P_1$ and $P_2$.
By Lemma~\ref{lemma2.1} and $P_1+P_2=Q$, $(P_1+y_1)\cap(P_2+y_2)=\emptyset$ for all $y_1\in Q$ and $y_2\in\mathbb{F}_q^{2\nu}\setminus Q$.
Therefore, the  desired result follows.
$\qed$

\begin{lemma}\label{lemma4.1.1}
Every $(1,0)$-spread in $\mathbb{F}_q^{2}$ is of type I.
\end{lemma}
\proof Let ${\cal S}$ be a $(1,0)$-spread in $\mathbb{F}_q^{2}$. If there exist $P_1+x_1,P_2+x_2\in{\cal S}$ such that  $P_1\not=P_2$,
from $(P_1+x_1)\cap(P_2+x_2)=\emptyset$ and Lemma~\ref{lemma2.1}, we deduce that $\dim((P_1+x_1)\vee(P_2+x_2))=3>2$, a contradiction. $\qed$

Now we give two useful lemmas by using a method in \cite{D'haeseleer}.

\begin{lemma}\label{lemma3.8}
 All the characteristic vectors of  $(\nu,0)$-spreads of type I in $\mathbb{F}_q^{2\nu}$ form a basis for the space  $V_{(0,0)}\perp V_{(1,0)}\perp\cdots\perp V_{(\nu,0)}$.
\end{lemma}
\proof Let ${\cal S}$ be a $(\nu,0)$-spread of type I in $\mathbb{F}_q^{2\nu}$. By Lemma~\ref{lemma3.1},
there exist $v_{(j,\eta)}\in V_{(j,\eta)}$ and $a_{(j,\eta)}\in\mathbb{R}$  such that the characteristic vector of ${\cal S}$
 $$\chi_{{\cal S}}=\sum_{j=0}^{\nu}a_{(j,0)}v_{(j,0)}+\sum_{j=0}^{\nu-1}a_{(j,1)}v_{(j,1)}.$$ Let $$u=(u_{(0,0)},u_{(0,1)},\ldots,u_{(\nu-1,0)},u_{(\nu-1,1)},u_{(\nu,0)})$$
 be the inner distribution of ${\cal S}$. Then $u_{(j,\eta)}=0$ for all $j\not=0$ and $\eta=0,1$, and
 $$u_{(0,\eta)}=\frac{|R_{(0,\eta)}\cap({\cal S}\times {\cal S})|}{|{\cal S}|}
 =\left\{\begin{array}{ll}
 1 & \hbox{if}\;\eta=0,\\
 q^{\nu}-1 & \hbox{if}\;\eta=1.
 \end{array}\right.$$
 Write $u(q_{(j,\eta)}(i,\xi))=(w_{(0,0)},w_{(0,1)},\ldots,w_{(\nu-1,0)},w_{(\nu-1,1)},w_{(\nu,0)})$, where $(q_{(j,\eta)}(i,\xi))$ is the  second eigenmatrix of ${\mathfrak X}({\cal O}_{\nu})$.
 By Lemmas~\ref{lemma3.3},~\ref{lemma3.4} and~(\ref{equa3.1}), we obtain
 \begin{eqnarray*}\label{equa3.4}
 w_{(j,\eta)}&=&q_{(j,\eta)}(0,0)+(q^{\nu}-1)q_{(j,\eta)}(0,1)
 =m_{(j,\eta)}(1+p_{(0,1)}(j,\eta))\\
 &=&\left\{\begin{array}{ll}
 q^{\nu}m_{(j,\eta)} & \hbox{if}\;\eta=0,\\
 0 & \hbox{if}\;\eta=1.
 \end{array}\right.
 \end{eqnarray*}
 By Lemma~\ref{lemma3.1} again,  $a_{(j,1)}=0$ for all $j\not=0$. It follows that $\chi_{{\cal S}}\in V_{(0,0)}\perp V_{(1,0)}\perp\cdots\perp V_{(\nu,0)}$.
 Since no two distinct $(\nu,0)$-spreads of type I in $\mathbb{F}_q^{2\nu}$ contain the same $(\nu,0)$-flat,
 we know that all the characteristic vectors of $(\nu,0)$-spreads of type I in $\mathbb{F}_q^{2\nu}$ are
linearly independent. Since $\dim(V_{(0,0)}\perp V_{(1,0)}\perp\cdots\perp V_{(\nu,0)})=\prod_{t=1}^{\nu}(q^{t+e-1}+1)$
is equal to the number of all $(\nu,0)$-spreads of type I in $\mathbb{F}_q^{2\nu}$, the desired result follows. $\qed$

\begin{lemma}\label{lemma3.9}
Let $\nu\geq 2$ and $(\nu,q,e)\not=(2,2,0)$.
Then the space $V_{(0,0)}\perp V_{(1,0)}\perp V_{(1,1)}\perp\cdots\perp V_{(\nu-1,0)}\perp V_{(\nu-1,1)}\perp V_{(\nu,0)}$
is spanned by all the characteristic vectors of $(\nu,0)$-spreads of type I\!I in $\mathbb{F}_q^{2\nu}$. Moreover,
 $E_{(j,\eta)}\chi_{{\cal S}}\not=0$ for all $j\not=0$ and $\eta=0,1$.
\end{lemma}
\proof Suppose that ${\cal S}$ is a $(\nu,0)$-spread of type I\!I in $\mathbb{F}_q^{2\nu}$.
 Let $$u=(u_{(0,0)},u_{(0,1)},\ldots,u_{(\nu-1,0)},u_{(\nu-1,1)},u_{(\nu,0)})$$
 be the inner distribution of ${\cal S}$. Then $u_{(1,0)}=0$ and $u_{(j,\eta)}=0$ for all $j\geq2$ and $\eta=0,1$, and
 \begin{eqnarray*}
 u_{(0,1)}&=&\frac{|R_{(0,1)}\cap({\cal S}\times {\cal S})|}{|{\cal S}|}
 =\frac{q^{2\nu}-2q^{\nu+1}-q^{\nu}+2q^2}{q^{\nu}}=q^\nu-2q-1+\frac{2q^2}{q^\nu},\\
 u_{(1,1)}&=&\frac{|R_{(1,1)}\cap({\cal S}\times {\cal S})|}{|{\cal S}|}
 =\frac{2q(q^{\nu}-q)}{q^{\nu}}=\frac{2q^2(q^{\nu-1}-1)}{q^\nu}.
 \end{eqnarray*}
 Write $u(q_{(j,\eta)}(i,\xi))=(w_{(0,0)},w_{(0,1)},\ldots,w_{(\nu-1,0)},w_{(\nu-1,1)},w_{(\nu,0)}).$
 By Lemmas~\ref{lemma3.3},~\ref{lemma3.4} and~(\ref{equa3.1}), we obtain
 \begin{eqnarray*}\label{equa3.4}
 w_{(j,\eta)}&=&q_{(j,\eta)}(0,0)+u_{(0,1)}q_{(j,\eta)}(0,1)+u_{(1,1)}q_{(j,\eta)}(1,1)\\
 &=&m_{(j,\eta)}\left(1+\frac{u_{(0,1)}p_{(0,1)}(j,\eta))}{v_{(0,1)}}+\frac{u_{(1,1)}p_{(1,1)}(j,\eta))}{v_{(1,1)}}\right)\\
 &=&\left\{\begin{array}{ll}
 m_{(j,0)}\left(1+u_{(0,1)}+\frac{u_{(1,1)}q(q^{\nu-1}-1)p_1^{(2\nu)}(j)}{v_{(1,1)}}\right) &\hbox{if}\;\eta=0,\\
 m_{(j,1)}\left(1-\frac{u_{(0,1)}}{q^\nu-1}-\frac{u_{(1,1)}qp_1^{(2(\nu-1))}(j)}{v_{(1,1)}}\right) &\hbox{if}\;\eta=1
 \end{array}\right.  \\
 &=&\left\{\begin{array}{ll}
 m_{(j,0)}\left(q^\nu-2q+\frac{2}{q^{\nu-2}}+\frac{2(q^{\nu-1}-1)\left(q^{e}{\nu-j\brack 1}_q-{j\brack 1}_q\right)}{q^{\nu-2+e}{\nu\brack 1}_q}\right) &\hbox{if}\;\eta=0,\\
 m_{(j,1)}\left(\frac{2q}{q^\nu-1}-\frac{2}{q^{\nu-2}(q^\nu-1)}-\frac{2\left(q^e{\nu-1-j\brack 1}_q-{j\brack 1}_q\right)}{q^{\nu-2+e}{\nu\brack 1}_q}\right) &\hbox{if}\;\eta=1.
 \end{array}\right.
 \end{eqnarray*}
 For all $0\leq j\leq \nu-1$, from $q^{e}{\nu-j\brack 1}_q-{j\brack 1}_q>-{\nu\brack 1}_q$, we deduce that
 $$\frac{w_{(j,0)}}{m_{(j,0)}}>\frac{w_{(\nu,0)}}{m_{(\nu,0)}}=\frac{(q^{\nu-1}-2)^2}{q^{\nu-2}}+\frac{2(q^{\nu-1}-1)}{q^{\nu-2}}\left(1-\frac{1}{q^e}\right)\geq0.
 $$
 Note that $\frac{w_{(\nu,0)}}{m_{(\nu,0)}}=0$ if and only if $(\nu,q,e)=(2,2,0)$.
 For all $1\leq j\leq \nu-1$, from $q^{e}{\nu-1-j\brack 1}_q-{j\brack 1}_q< q^e{\nu-2\brack 1}_q$, we deduce that
 $$
 \frac{w_{(j,1)}}{m_{(j,1)}}>\frac{2q}{q^\nu-1}-\frac{2}{q^{\nu-2}(q^\nu-1)}-\frac{2(q^{\nu-2}-1)}{q^{\nu-2}(q^\nu-1)}
 =\frac{2(q-1)}{q^\nu-1}>0.
 $$
 Note that $ \frac{w_{(0,1)}}{m_{(0,1)}}=0$.
 By Lemma~\ref{lemma3.1}, we have $$\chi_{\cal S}\in V_{(0.0)}\perp V_{(1,0)}\perp V_{(1,1)}\perp\cdots\perp V_{(\nu-1,0)}\perp V_{(\nu-1,1)}\perp V_{(\nu,0)}$$ and $E_{(j,\eta)}\chi_{{\cal S}}\not=0$ for all $j\not=0$ and $\eta=0,1$.

 Next, we show that the space $V_{(0.0)}\perp V_{(1,0)}\perp V_{(1,1)}\perp\cdots\perp V_{(\nu-1,0)}\perp V_{(\nu-1,1)}\perp V_{(\nu,0)}$
is spanned by all the characteristic vectors of $(\nu,0)$-spreads of type I\!I in $\mathbb{F}_q^{2\nu}$.
Suppose that these characteristic vectors would span $V_{(0,0)}\perp W$ with $V_{(1,0)}\perp V_{(1,1)}\perp\cdots\perp V_{(\nu-1,0)}\perp V_{(\nu-1,1)}\perp V_{(\nu,0)}=W\perp U$.
If there exists a $\psi\in U\setminus\{0\}$, then there exists some $(j,\eta)$ with $j\not=0$ and $\eta=0,1$, such that $E_{(j,\eta)}\psi\not=0$.
Suppose that ${\cal S}$ is a $(\nu,0)$-spread of type I\!I in $\mathbb{F}_q^{2\nu}$.
Then $\chi_{\cal S}\in V_{(0.0)}\perp W$. It follows that $\chi_{\cal S}\cdot\psi=0$. By the transitivity of $AG_{2\nu}$,   elements of $AG_{2\nu}$ carry  $(\nu,0)$-spreads of  type I\!I in $\mathbb{F}_q^{2\nu}$ to $(\nu,0)$-spreads of
type  I\!I. Therefore, we  obtain
$$\psi\cdot\chi_{{\cal S}^g}=0\quad\hbox{for all}\; g\in AG_{2\nu}.$$
By Lemmas~\ref{lemma3.1.2} and~\ref{lemma-trs},  $E_{(j,\eta)}\chi_{\cal S}=0$ for such $(j,\eta)$, a contradiction.$\qed$

Next, we give several equivalent definitions for a Cameron-Liebler set in ${\cal O}_{\nu}$.

\begin{thm}\label{Thm-eql}
Let ${\cal L}$ be a non-empty set in ${\cal O}_{\nu}$
with $|{\cal L}|=x\prod_{i=1}^{\nu}(q^{i+e-1}+1)$. Then the following properties are equivalent.
\begin{itemize}
\item[\rm(i)]
$\chi_{{\cal L}}\in {\rm Im}(M^{\top})$.

\item[\rm(ii)]
$\chi_{{\cal L}}\in \ker(M)^{\perp}$.

\item[\rm(iii)]
$\chi_{{\cal L}}\in V_{(0,0)}\perp V_{(0,1)}$, where $V_{(0,0)}=\langle{\rm j}\rangle$.

\item[\rm(iv)]
The vector $v=\chi_{{\cal L}}-xq^{-\nu}{\rm j}\in V_{(0,1)}$.

\item[\rm(v)]
For every $F\in{\cal O}_{\nu}$, the number of elements $F'$ in ${\cal L}$ satisfying $(F,F')\in R_{(1,\xi)}$  is
$$\left\{\begin{array}{ll}
q^{1+e}{\nu-1\brack 1}_q+xq^{e}  &\hbox{if}\;F\in{\cal L}\;\hbox{and}\;\xi=0,\\
xq^{e}    &\hbox{if}\;F\not\in{\cal L}\;\hbox{and}\;\xi=0,\\
(x-1)q^{1+e}{\nu-1\brack 1}_q  &\hbox{if}\;F\in{\cal L}\;\hbox{and}\;\xi=1,\\
xq^{1+e}{\nu-1\brack 1}_q  &\hbox{if}\;F\not\in{\cal L}\;\hbox{and}\;\xi=1.
\end{array}\right.
$$
\end{itemize}

If $(\nu,q,e)\not=(2,2,0)$, then the following properties are equivalent to the previous ones.
\begin{itemize}
\item[\rm(vi)]
For every $(\nu,0)$-conjugate switching set ${\cal R}$ and ${\cal R}'$ in $\mathbb{F}_q^{2\nu}$, we have $|{\cal L}\cap{\cal R}|$ = $|{\cal L}\cap{\cal R}^{\prime}|$.

\item[\rm(vii)]
  $|{\cal L}\cap{\cal S}| = x$ for every $(\nu,0)$-spread ${\cal S}$ in $\mathbb{F}_q^{2\nu}$.
\end{itemize}
\end{thm}
\proof (i) $\Leftrightarrow$ (ii)  Trivial.

(ii) $\Leftrightarrow$ (iii)
By Lemma~\ref{lemma2.5}, the matrix $M$ has $q^{2\nu}$
rows with rank$(M)=(q^{\nu}-1)(q^{\nu+e-1}+1)+1$. Since
$\dim V_{(0,1)}=(q^{\nu}-1)(q^{\nu+e-1}+1)$, we have $\dim(V_{(0,0)}\perp V_{(0,1)})=\dim{\rm Im}(M^{\top})$.
For any $\{a\}\in {\cal O}_0$, we show that $\chi_{{\cal O}'_{\nu}(\{a\})}\in V_{(0,0)}\perp V_{(0,1)}$, where  ${\cal O}_{\nu}'(\{a\})$ is the set
of all maximal totally isotropic flats in  ${\cal O}_\nu$ containing the fixed $(0,0)$-flat $\{a\}$.
By Lemma~\ref{lemma3.1},
there exist $v_{(j,\eta)}\in V_{(j,\eta)}$ and $a_{(j,\eta)}\in\mathbb{R}$  such that
 $$\chi_{{\cal O}'_{\nu}(\{a\})}=\sum_{j=0}^{\nu}a_{(j,0)}v_{(j,0)}+\sum_{j=0}^{\nu-1}a_{(j,1)}v_{(j,1)}.$$
 In order to prove $\chi_{{\cal O}'_{\nu}(\{a\})}\in V_{(0,0)}\perp V_{(0,1)}$, we only need to show that $a_{(j,\eta)}=0$ for all $j\not=0$ and $\eta=0,1$.

Let $$u=(u_{(0,0)},u_{(0,1)},\ldots,u_{(\nu-1,0)},u_{(\nu-1,1)},u_{(\nu,0)})$$
 be the inner distribution of ${\cal O}'_{\nu}(\{a\})$. By Lemma~\ref{lemma3.2}, the valency $v^{(2\nu)}_i$ of ${\mathfrak X}({\cal M}_{\nu})$
  is $q^{i(i+2e-1)/2}{\nu\brack i}_q$. By Lemmas~\ref{lemma2.1} and~\ref{lemma2.4}, we have $u_{(i,1)}=0$ for all $i=0,1,\ldots,\nu-1$, and
 $$u_{(i,0)}=\frac{|R_{(i,0)}\cap({\cal O}'_{\nu}(\{a\})\times {\cal O}'_{\nu}(\{a\}))|}{|{\cal O}'_{\nu}(\{a\})|}
 =v^{(2\nu)}_i=q^{i(i+2e-1)/2}{\nu\brack i}_q\;\hbox{for all}\;i=0,1,\ldots,\nu.$$
 Write $u(q_{(j,\eta)}(i,\xi))=(w_{(0,0)},w_{(0,1)},\ldots,w_{(\nu-1,0)},w_{(\nu-1,1)},w_{(\nu,0)})$, where $(q_{(j,\eta)}(i,\xi))$ is the  second eigenmatrix of ${\mathfrak X}({\cal O}_{\nu})$.
 By Lemma~\ref{lemma3.3} and~(\ref{equa3.1}), we obtain
 \begin{equation}\label{equa3.4}
 w_{(j,\eta)}=\sum_{i=0}^{\nu}v^{(2\nu)}_iq_{(j,\eta)}(i,0)
 =\sum_{i=0}^{\nu}\frac{v^{(2\nu)}_im_{(j,\eta)}p_{(i,0)}(j,\eta)}{v_{(i,0)}}
 =\sum_{i=0}^{\nu}\frac{m_{(j,\eta)}p_{(i,0)}(j,\eta)}{q^i}.
 \end{equation}

{\it Case}~1: $\eta=0$. By Lemma~\ref{lemma3.4}, (\ref{equa3.1.1}) and (\ref{equa3.4}), we have
 $$
 w_{(j,0)}=\sum_{i=0}^{\nu}\frac{m_{(j,0)}p_{(i,0)}(j,0)}{q^i}
 =m^{(2\nu)}_{j}\sum_{i=0}^{\nu}p_{i}^{(2\nu)}(j)
 =\left\{\begin{array}{ll}
 \prod_{t=1}^{\nu}(q^{t+e-1}+1) & \hbox{if}\; j=0,\\
 0   & \hbox{otherwise}.
 \end{array}
 \right.
 $$

 {\it Case}~2: $\eta=1$. By Lemma~\ref{lemma3.4}, (\ref{equa3.1.1}) and (\ref{equa3.4}), we have
 \begin{eqnarray*}
 w_{(j,1)}&=&\sum_{i=0}^{\nu}\frac{m_{(j,1)}p_{(i,0)}(j,1)}{q^i}\\
 &=&(q^{\nu}-1)(q^{\nu+e-1}+1)m_j^{(2\nu-2)}\sum_{i=0}^{\nu-1}p_{i}^{(2\nu-2)}(j)\\
 &=&\left\{\begin{array}{ll}
 (q^{\nu}-1)\prod_{t=1}^{\nu}(q^{t+e-1}+1) & \hbox{if}\; j=0,\\
 0   & \hbox{otherwise}.
 \end{array}
 \right.
 \end{eqnarray*}

By Lemma~\ref{lemma3.1}, we obtain $w_{(j,\eta)}=0$ for all $j\not=0$ and $\eta=0,1$,  which implies that $a_{(j,\eta)}=0$ for all $j\not=0$ and $\eta=0,1$. It follows that $\chi_{{\cal O}'_{\nu}(\{a\})}\in V_{(0,0)}\perp V_{(0,1)}$. From ${\rm Im}(M^{\top})\subseteq V_{(0,0)}\perp V_{(0,1)}$ and $\dim(V_{(0,0)}\perp V_{(0,1)})=\dim{\rm Im}(M^{\top})$, we deduce that ${\rm Im}(M^{\top})=V_{(0,0)}\perp V_{(0,1)}$.

(iii) $\Leftrightarrow$ (iv) From $\chi_{{\cal L}}\in V_{(0,0)}\perp V_{(0,1)}$, we deduce that $v=\chi_{{\cal L}}-xq^{-\nu}{\rm j}\in V_{(0,0)}\perp V_{(0,1)}$. By Lemma~\ref{lemma2.4} (i),
we have
 $$v\cdot {\rm j} = |{\cal L}| - xq^{-\nu}|{\cal O}_{\nu}|=x\prod_{i=1}^{\nu}(q^{i+e-1}+1) - xq^{-\nu}q^{\nu}\prod_{i=1}^{\nu}(q^{i+e-1}+1) = 0,$$
 which implies that $v\in V_{(0,1)}$ by $V_{(0,0)}=\langle{\rm j}\rangle$. Conversely, if $v=\chi_{{\cal L}}-xq^{-\nu}{\rm j}\in V_{(0,1)}$,
 then $\chi_{{\cal L}}\in V_{(0,0)}\perp V_{(0,1)}$ is obvious.

 (iv) $\Leftrightarrow$ (v) Note that the matrix $A_{(1,\xi)}$ corresponds to the relation $R_{(1,\xi)}$,
 and the eigenvalue of $A_{(1,\xi)}$ corresponding to $V_{(j,\eta)}$ is
$p_{(1,\xi)}(j,\eta)$.
  So, $(A_{(1,\xi)}\chi_{{\cal L}})_{F}$ gives the
number of elements $F'$ in ${\cal L}$ satisfying $(F,F')\in R_{(1,\xi)}$. Since $v=\chi_{{\cal L}}-xq^{-\nu}{\rm j}\in V_{(0,1)}$,
$$A_{(1,\xi)}\chi_{{\cal L}} = A_{(1,\xi)}v+xq^{-\nu}A_{(1,\xi)}{\rm j}=p_{(1,\xi)}(0,1)v+xq^{-\nu}p_{(1,\xi)}(0,0){\rm j}.$$

{\it Case}~1: $\xi=0$. By Lemma~\ref{lemma3.4}, we have
\begin{eqnarray*}
A_{(1,0)}\chi_{{\cal L}}  &=&q^{1+e}{\nu-1\brack 1}_q(\chi_{{\cal L}}-xq^{-\nu}{\rm j})   +xq^{-\nu}q^{1+e}{\nu\brack 1}_q{\rm j}\\
  &=&q^{1+e}{\nu-1\brack 1}_q\chi_{{\cal L}}+xq^{e}{\rm j}.
\end{eqnarray*}

{\it Case}~2: $\xi=1$. By Lemma~\ref{lemma3.4}, we have
\begin{eqnarray*}
A_{(1,1)}\chi_{{\cal L}}  &=&-q^{1+e}{\nu-1\brack 1}_q(\chi_{{\cal L}}-xq^{-\nu}{\rm j})   +xq^{-\nu}q^{1+e}{\nu\brack 1}_q(q^{\nu-1}-1){\rm j}\\
  &=&-q^{1+e}{\nu-1\brack 1}_q\chi_{{\cal L}}+xq^{1+e}{\nu-1\brack 1}_q{\rm j}\\
  &=&q^{1+e}{\nu-1\brack 1}_q(x{\rm j}-\chi_{{\cal L}}).
\end{eqnarray*}

Conversely, from property (v), we deduce that
$$
A_{(1,\xi)}\chi_{{\cal L}}=\left\{\begin{array}{ll}
q^{1+e}{\nu-1\brack 1}_q\chi_{{\cal L}}+xq^{e}{\rm j}  &\hbox{if}\;\xi=0,\\
-q^{1+e}{\nu-1\brack 1}_q\chi_{{\cal L}}+xq^{1+e}{\nu-1\brack 1}_q{\rm j}  &\hbox{if}\;\xi=1.
\end{array}\right.
$$
If $\xi=0$, by Lemmas~\ref{lemma3.3} and~\ref{lemma3.4}, we have
\begin{eqnarray*}
&&A_{(1,\xi)}(\chi_{{\cal L}}-xq^{-\nu}{\rm j})\\
&=&q^{1+e}{\nu-1\brack 1}_q\chi_{{\cal L}}+x\left(q^{e}
-q^{-\nu}q^{1+e}{\nu\brack 1}_q\right){\rm j} \\
&=&q^{1+e}{\nu-1\brack 1}_q(\chi_{{\cal L}}-xq^{-\nu}{\rm j})
=p_{(1,0)}(0,1)(\chi_{{\cal L}}-xq^{-\nu}{\rm j}).
\end{eqnarray*}
If $\xi=1$, by Lemmas~\ref{lemma3.3} and~\ref{lemma3.4}, we have
\begin{eqnarray*}
&&A_{(1,\xi)}(\chi_{{\cal L}}-xq^{-\nu}{\rm j})\\
&=&-q^{1+e}{\nu-1\brack 1}_q\chi_{{\cal L}}+x\left(q^{1+e}{\nu-1\brack 1}_q
-q^{-\nu}(q^{\nu-1}-1)q^{1+e}{\nu\brack 1}_q\right){\rm j} \\
&=&-q^{1+e}{\nu-1\brack 1}_q(\chi_{{\cal L}}-xq^{-\nu}{\rm j})
=p_{(1,1)}(0,1)(\chi_{{\cal L}}-xq^{-\nu}{\rm j}).
\end{eqnarray*}
It follows that
$\chi_{{\cal L}}-xq^{-\nu}{\rm j}$
is an eigenvector for $A_{(1,\xi)}$ with eigenvalue $p_{(1,\xi)}(0,1)$.
By Lemma~\ref{lemma3.9-N}, $\chi_{{\cal L}}-xq^{-\nu}{\rm j}\in V_{(0,1)}$.

To end this proof, we show that the properties (vi) and (vii) are equivalent with the other properties.

(ii) $\Rightarrow$ (vi) Since the $(\nu,0)$-conjugate switching set ${\cal R}$ and ${\cal R}'$  in $\mathbb{F}_q^{2\nu}$ covers the same subset of $\mathbb{F}_q^{2\nu}$, we have $\chi_{{\cal R}} - \chi_{{\cal R}'}\in \ker(M)$, which implies that
$$0 = \chi_{{\cal L}}\cdot(\chi_{{\cal R}} - \chi_{{\cal R}'}) = \chi_{{\cal L}}\cdot\chi_{{\cal R}} - \chi_{{\cal L}}\cdot\chi_{{\cal R}'} = |{\cal L}\cap {\cal R}| -  |{\cal L}\cap {\cal R}'|,$$ as desired.

(vi) $\Rightarrow$ (vii) For any two $(\nu,0)$-spreads ${\cal S}$ and ${\cal S}'$ in  $\mathbb{F}_q^{2\nu}$, both
${\cal S}\setminus{\cal S}'$ and ${\cal S}'\setminus{\cal S}$ form a $(\nu,0)$-conjugate switching set  in  $\mathbb{F}_q^{2\nu}$.
So $|{\cal L}\cap ({\cal S}\setminus{\cal S}')|=|{\cal L}\cap ({\cal S}'\setminus{\cal S})|$,
which implies that $|{\cal L}\cap {\cal S}|=|{\cal L}\cap{\cal S}')|= c$.

Now we prove $c = x = |{\cal L}|/\prod_{i=1}^{\nu}(q^{i+e-1}+1)$.
Note that $\{P+y: y\in \mathbb{F}_q^{2\nu}\}$ is a $(\nu,0)$-spread in $\mathbb{F}_q^{2\nu}$ for every $P\in{\cal M}_{\nu}$,
where ${\cal M}_{\nu}$ is the set of all maximal totally isotropic subspaces in $\mathbb{F}_q^{2\nu}$. Since ${\cal L}\subseteq\bigcup_{P\in{\cal M}_{\nu}}\{P+y:y\in \mathbb{F}_q^{2\nu}\}$, by Lemma~\ref{lemma2.4} (i), we have
\begin{eqnarray*}
  x\prod_{i=1}^{\nu}(q^{i+e-1}+1)&=&|{\cal L}|=\left|{\cal L}\cap\bigcup_{P\in{\cal M}_{\nu}}\{P+y:y\in \mathbb{F}_q^{2\nu}\}\right|\\
 &=& c|{\cal M}_{\nu}|=c\prod_{i=1}^{\nu}(q^{i+e-1}+1),
\end{eqnarray*}
which implies that $x = c$.

(vii) $\Rightarrow$ (iii) Let ${\cal S}$ be a $(\nu,0)$-spread  in $\mathbb{F}_q^{2\nu}$. Since $|{\cal L}\cap{\cal S}| = x$, we have $\chi_{\cal L}\cdot\chi_{\cal S}=x$.
By the transitivity of $AG_{2\nu}$,   elements of $AG_{2\nu}$ carry $(\nu,0)$-spreads of type I and type I\!I in $\mathbb{F}_q^{2\nu}$ to $(\nu,0)$-spreads of
type I and type I\!I,  respectively. Therefore, we obtain
$$\chi_{\cal L}\cdot\chi_{{\cal S}^g}=x\quad\hbox{for all}\; g\in AG_{2\nu}.$$ There are the following two cases to be considered.

{\it Case}~1: $\nu=1$. By Lemma~\ref{lemma4.1.1}, ${\cal S}$ is a $(1,0)$-spread of type I.
From Lemma~\ref{lemma3.8}, we deduce that $E_{(1,0)}\chi_{\cal S}\not=0$, which implies that $E_{(1,0)}\chi_{\cal L}=0$ by
Lemmas~\ref{lemma3.1.2} and~\ref{lemma-trs}. By Lemma~\ref{lemma3.1}, $\chi_{\cal L}\in V_{(0,0)}\perp V_{(0,1)}$.

{\it Case}~2: $\nu\geq2$. Suppose that ${\cal S}$ is a $(\nu,0)$-spread of type I\!I. From Lemma~\ref{lemma3.9},
we deduce that $E_{(j,\eta)}\chi_{\cal S}\not=0$ for all $j\not=0$ and $\eta=0,1$. By
Lemmas~\ref{lemma3.1.2} and~\ref{lemma-trs},  $E_{(j,\eta)}\chi_{\cal L}=0$  for all $j\not=0$ and $\eta=0,1$. So, $\chi_{\cal L}\in V_{(0,0)}\perp V_{(0,1)}$ by Lemma~\ref{lemma3.1}.
$\qed$

The following result is a useful property of Cameron-Liebler sets in ${\cal O}_{\nu}$.

\begin{lemma}\label{lemma3.10}
Let ${\cal L}$ be a Cameron-Liebler set with parameter $x$ in ${\cal O}_{\nu}$.
For every $F\in{\cal O}_{\nu}$, the number of elements $F'$ in ${\cal L}$ satisfying $(F,F')\in R_{(i,\xi)}$  is
\begin{equation}\label{equa8}
\left\{\begin{array}{ll}
q^{i(i+2e+1)/2}{\nu-1\brack i}_q+xq^{i(i+2e-1)/2}{\nu-1\brack i-1}_q  &\hbox{if}\;F\in{\cal L}\;\hbox{and}\;\xi=0,\\
xq^{i(i+2e-1)/2}{\nu-1\brack i-1}_q   &\hbox{if}\;F\not\in{\cal L}\;\hbox{and}\;\xi=0,\\
(x-1)q^{i(i+2e+1)/2}{\nu-1\brack i}_q  &\hbox{if}\;F\in{\cal L}\;\hbox{and}\;\xi=1,\\
xq^{i(i+2e+1)/2}{\nu-1\brack i}_q  &\hbox{if}\;F\not\in{\cal L}\;\hbox{and}\;\xi=1.
\end{array}\right.
\end{equation}
\end{lemma}
\proof The proof  is similar to that of Theorem~\ref{Thm-eql}, and will be omitted. $\qed$

\medskip
\noindent{\bf Remarks.} Let ${\cal L}$ be a non-empty set in ${\cal O}_{\nu}$ with $|{\cal L}|=x\prod_{i=1}^{\nu}(q^{i+e-1}+1)$.
Suppose that $(i,\xi)\not=(0,0)$, and the number of elements $F'$ in ${\cal L}$ satisfying $(F,F')\in R_{(i,\xi)}$  is
as in (\ref{equa8}) for every $F\in{\cal O}_{\nu}$. If the eigenvalue $p_{(i,\xi)}(0,1)$ of $A_{(i,\xi)}$ only corresponds to $V_{(0,1)}$,
then ${\cal L}$ is a Cameron-Liebler set with parameter $x$ in ${\cal O}_{\nu}$.
By Lemmas~\ref{lemma3.9-N} and~\ref{lemma3.10},  we may obtain several equivalent definitions for a Cameron-Liebler set in ${\cal O}_{\nu}$, 
except in the following cases: {\rm(a)} $\nu\geq2$ and $(i,\xi)=(0,1)$; {\rm(b)} $\nu\geq2$ and $(i,\xi)=(\nu,0)$; 
{\rm(c)}  $i$ is even and $2\leq i\leq\nu-1$.

\section{Classification results}
In this section, we give some classification results for Cameron-Liebler sets in ${\cal O}_{\nu}$.
We begin with some properties of Cameron-Liebler sets in ${\cal O}_{\nu}$.

By Theorem~\ref{Thm-eql}, we obtain the following result.

\begin{lemma}\label{lemma5.1}
Let ${\cal L}$ and ${\cal L}'$ be two Cameron-Liebler sets in ${\cal O}_{\nu}$ with parameters $x$ and $x'$, respectively. Then the following hold.
\begin{itemize}
\item[\rm(i)]
$0\leq x\leq q^{\nu}$.

\item[\rm(ii)]
The set of all flats in ${\cal O}_{\nu}$ not in ${\cal L}$ is a Cameron-Liebler set in  ${\cal O}_{\nu}$ with parameter $q^{\nu}-x$.

\item[\rm(iii)]
If ${\cal L}\cap{\cal L}'=\emptyset$, then ${\cal L}\cup{\cal L}'$ is a Cameron-Liebler set in  ${\cal O}_{\nu}$ with parameter $x+x'$.

\item[\rm(iv)]
If ${\cal L}'\subseteq{\cal L}$, then ${\cal L}\setminus{\cal L}'$ is a Cameron-Liebler set in  ${\cal O}_{\nu}$  with parameter $x-x'$.
\end{itemize}
\end{lemma}

Recall that  ${\cal O}_{\nu}'(\{a\})$ is the set
of all maximal totally isotropic flats in  ${\cal O}_\nu$ containing a fixed $(0,0)$-flat $\{a\}$.
It is well known that  the set ${\cal O}'_{\nu}(\{a\})$ is a Cameron-Liebler set in  ${\cal O}_{\nu}$ with parameter $1$,
and the set ${\cal O}_{\nu}$ is a Cameron-Liebler set in  ${\cal O}_{\nu}$ with parameter $q^\nu$.

If $\nu=1$, we have the following classification result.

 \begin{thm}\label{thm5.2}
 Let ${\cal P}$ be the collection of all  $(1,0)$-spreads in $\mathbb{F}_q^{2}$, and ${\cal L}$ be a Cameron-Liebler set in ${\cal O}_{1}$. Then
 the following statements are equivalent.
 \begin{itemize}
 \item[\rm(i)]
 ${\cal L}$ is a Cameron-Liebler set in ${\cal O}_{1}$ with parameter $x$.

 \item[\rm(ii)]
 $|{\cal L}\cap{\cal S}|=x$ for all ${\cal S}\in{\cal P}$ and ${\cal L}=\bigcup_{{\cal S}\in{\cal P}}({\cal L}\cap{\cal S})$.
 \end{itemize}
 \end{thm}
 \proof By Lemma~\ref{lemma4.1.1}, every $(1,0)$-spread in $\mathbb{F}_q^{2}$ is of type I.
Suppose that $\{P_1,\ldots,P_{q^e+1}\}$ is the set of all subspaces of type $(1,0)$ in $\mathbb{F}_q^2$. Then ${\cal P}=\{{\cal S}_1,\ldots,{\cal S}_{q^e+1}\}$, where
 ${\cal S}_i=\{P_i+x_i : x_i\in\mathbb{F}_q^2\}$ for all $i=1,\ldots,q^e+1$.

 (i) $\Rightarrow$ (ii) By Theorem~\ref{Thm-eql}, we have $|{\cal L}\cap{\cal S}_i|=x$ for all $i=1,\ldots,q^e+1$.
 Since ${\cal S}_i\cap{\cal S}_j=\emptyset$
 for all $i\not=j$, the size of $\bigcup_{i=1}^{q^e+1}({\cal L}\cap{\cal S}_i)$ is $x(q^e+1)$. From ${\cal L}\supseteq\bigcup_{i=1}^{q^e+1}({\cal L}\cap{\cal S}_i)$
 and $|{\cal L}|=x(q^e+1)$, we deduce that ${\cal L}=\bigcup_{i=1}^{q^e+1}({\cal L}\cap{\cal S}_i)$.

 (ii) $\Rightarrow$ (i)
 Since $|{\cal L}\cap{\cal S}_i|=x$ for all $i=1,\ldots,q^e+1$, and
 ${\cal S}_i\cap{\cal S}_j=\emptyset$
 for all $i\not=j$, we obtain that the size of $\bigcup_{i=1}^{q^e+1}({\cal L}\cap{\cal S}_i)$ is $x(q^e+1)$. By Theorem~\ref{Thm-eql},
 ${\cal L}$ is a Cameron-Liebler set in ${\cal O}_{1}$ with parameter $x$. $\qed$

Next, we consider Cameron-Liebler sets in ${\cal O}_{\nu}$ for $\nu\geq2$.
Firstly, we give a classification result for Cameron-Liebler sets in ${\cal O}_{\nu}$ with parameter $x=1$.
We begin with two useful lemmas.

\begin{lemma}\label{lemma5.3}{\rm(See \cite{Cameron2}).}
Let $C$ be a clique and $A$ a coclique in a vertex-transitive graph on $v$
vertices. Then $|C||A|\leq v$. Equality implies that $|C\cap A|=1$.
\end{lemma}

A family ${\cal F}\subseteq{\cal O}_{\nu}$ is called {\it intersecting} if $F\cap F'\not=\emptyset$
for all $F,F'\in{\cal F}$.

\begin{lemma}\label{lemma5.4}
Let ${\cal F}\subseteq{\cal O}_{\nu}$ be an intersecting family. Then $|{\cal F}|\leq\prod_{i=1}^{\nu}(q^{i+e-1}+1)$.
\end{lemma}
\proof We define a graph, denoted by $\Gamma$, on vertex set ${\cal O}_{\nu}$
by joining $F$ and $F'$ if they are intersecting. Since $AG_{2\nu}$ is a transitive permutation
group on ${\cal O}_{\nu}$, this graph is vertex-transitive. Let $\alpha(\Gamma)$ be the size of the largest
independent set in $\Gamma$.

Given a maximal totally isotropic subspace $P$ in $\mathbb{F}_q^{2\nu}$. Let ${\cal S}=\{P+x : x\in\mathbb{F}_q^{2\nu}\}$.
Then ${\cal S}$ is an independent set in $\Gamma$. So $\alpha(\Gamma)\geq q^{\nu}$.
Conversely, suppose $\alpha(\Gamma)> q^{\nu}$ and ${\cal I} = \{S_1, S_2, \ldots, S_{\alpha(\Gamma)}\}$ be a largest independent
set in $\Gamma$. Then $S_i\cap S_j=\emptyset$ for all $i\not=j$. Since $\bigcup_{i=1}^{\alpha(\Gamma)}S_i\subseteq\mathbb{F}_q^{2\nu}$ and
$\alpha(\Gamma)q^{\nu}>q^{2\nu}$, a contradiction. Hence $\alpha(\Gamma)=q^{\nu}$.
Since $\alpha(\Gamma)=q^{\nu}$ and $|{\cal O}_{\nu}|=q^{\nu}\prod_{i=1}^{\nu}(q^{i+e-1}+1)$, by Lemma~\ref{lemma5.3}, $|{\cal F}|\leq\prod_{i=1}^{\nu}(q^{i+e-1}+1)$. $\qed$

An intersecting family ${\cal F}\subseteq{\cal O}_{\nu}$ is called {\it maximum} if  $|{\cal F}|=\prod_{i=1}^{\nu}(q^{i+e-1}+1)$.

\begin{thm}\label{thm5.5}
Let ${\cal L}$ be a Cameron-Liebler set in ${\cal O}_{\nu}$ with parameter $x=1$. Then ${\cal L}$
is a maximum intersecting family.
\end{thm}
\proof For every $F\in{\cal L}$, by Lemma~\ref{lemma3.10}, the number of elements $F'$ in ${\cal L}$ satisfying $(F,F')\in \bigcup_{i=0}^{\nu-1}R_{(i,1)}$ is zero,
which implies that ${\cal L}$ is an intersecting family.  From $|{\cal F}|=\prod_{i=1}^{\nu}(q^{i+e-1}+1)$, we deduce that the desired result follows. $\qed$

Let $\nu\geq2$. We now discuss Cameron-Liebler sets in ${\cal O}_{\nu}$ with parameter $x\geq 2$. We begin with some useful notation and terminology.

From now on, we always assume that $1\leq i<\nu$ and $F$ is a fixed $(\nu+i,2i)$-flat in $ACG(2\nu,\,\mathbb{F}_q)$.
Let ${\cal O}_{\nu}(F)$ be the set of all $(\nu,0)$-flats contained in  $F$. Let ${\cal F}$ be the incidence matrix with rows indexed with vectors in $F$
and columns indexed with flats in ${\cal O}_{\nu}(F)$ such that the entry ${\cal F}_{x,S}=1$ if
and only if $x\in S$. The {\it characteristic vector} $f_{{\cal C}}$ of a subset ${\cal C}$ of ${\cal O}_{\nu}(F)$ is the column vector
in which its positions correspond to the elements of ${\cal O}_{\nu}(F)$, such that $(f_{{\cal C}})_{S}=1$ if $S\in{\cal O}_{\nu}(F)$ and 0 otherwise.
A subset ${\cal C}$ of ${\cal O}_{\nu}(F)$ is called a {\it Cameron-Liebler set} with parameter $x_F=|{\cal C}|/\prod_{j=1}^{i}(q^{j+e-1}+1)$
in ${\cal O}_{\nu}(F)$ if $f_{{\cal C}}\in {\rm Im}({\cal F}^{\top})$.

\begin{thm}\label{thm5.6}
Let ${\cal L}$ be a Cameron-Liebler set in ${\cal O}_{\nu}$ with parameter $x$. Then there exists an integer $x_F$ such that ${\cal L}\cap{\cal O}_{\nu}(F)$
is  a Cameron-Liebler set in ${\cal O}_{\nu}(F)$ with parameter $x_F$, where $0\leq x_F\leq\min\{x,q^i\}$. In particular,
${\cal O}_{\nu}(F)$ is  a Cameron-Liebler set in ${\cal O}_{\nu}(F)$ with parameter $x_F=q^i$.
\end{thm}
\proof Recall that $M$ is the incidence matrix with rows indexed with vectors in $\mathbb{F}_q^{2\nu}$
and columns indexed with flats in ${\cal O}_{\nu}$. After a suitable reordering of the
vectors and $(\nu,0)$-flats, we have
$$M=\left(\begin{array}{cc}
{\cal F} & M_2\\
0 & M_3\end{array}\right).$$
Since ${\cal L}$ is a Cameron-Liebler set in ${\cal O}_{\nu}$,
the  characteristic vector $\chi_{{\cal L}}\in{\rm Im}(M^{\top})$. Therefore,  there exists a vector
$\left(\begin{array}{c}
v_1\\ v_2\end{array}\right)$
such that $$\chi_{{\cal L}}=M^{\top}\left(\begin{array}{c}
v_1\\ v_2\end{array}\right)
=\left(\begin{array}{cc}
{\cal F}^{\top} & 0\\
M_2^{\top} & M_3^{\top}\end{array}\right)
\left(\begin{array}{c}
v_1\\ v_2\end{array}\right).$$
  It follows that
  $$f_{{\cal L}\cap{\cal O}_\nu(F)} = {\cal F}^{\top}v_1.$$ Therefore,
 ${\cal L}\cap{\cal O}_{\nu}(F)$ is  a Cameron-Liebler set in ${\cal O}_{\nu}(F)$.

Recall that a $(\nu,0)$-spread in $F$ is a set of pairwise disjoint $(\nu,0)$-flats in $F$ that partitions the set of all vectors in $F$.
Let ${\cal S}_1$ and ${\cal S}_2$ be two distinct $(\nu,0)$-spreads in $F$. Then both ${\cal S}_1\setminus{\cal S}_2$ and ${\cal S}_2\setminus{\cal S}_1$
is a $(\nu,0)$-conjugated switching set in $\mathbb{F}_q^{2\nu}$. By Theorem~\ref{Thm-eql}, we obtain
$$|{\cal L}\cap({\cal S}_1\setminus{\cal S}_2)|=|{\cal L}\cap({\cal S}_2\setminus{\cal S}_1)|,$$
which implies that $$|{\cal L}\cap{\cal O}_{\nu}(F)\cap{\cal S}_1|=|{\cal L}\cap{\cal S}_1|=|{\cal L}\cap{\cal S}_2|=|{\cal L}\cap{\cal O}_{\nu}(F)\cap{\cal S}_2|.$$
Therefore, there exists an integer $x_F$ such that $x_F=|{\cal L}\cap{\cal O}_{\nu}(F)\cap{\cal S}|$ for every $(\nu,0)$-spread ${\cal S}$ in $F$.

Without loss of generality, by the transitivity of $AG_{2\nu}$, we may assume that $F=\langle e_1,e_2,\ldots,e_{\nu+i}\rangle$,
where $e_j\,(1\leq j\leq\nu+i)$ is the row vector in $\mathbb{F}_q^{2\nu}$ in which its
$j$th coordinate is 1 and all other coordinates are 0. By \cite{Wan}, the number of maximal totally isotropic subspaces  contained in $F$ is $\varpi=\prod_{j=1}^{i}(q^{j+e-1}+1)$.
Let $\{P_1,P_2,\ldots,P_{\varpi}\}$ be the set of all maximal totally isotropic subspaces contained in $F$.
Then $\{P_s+y_s: y_s\in F\}\;(1\leq s\leq \varpi)$ are $(\nu,0)$-spreads  in $F$, and
${\cal O}_{\nu}(F)=\bigcup_{s=1}^{\varpi}\{P_s+y_s: y_s\in F\}$.
Since
$$
|{\cal L}\cap{\cal O}_{\nu}(F)|=\left|{\cal L}\cap\bigcup_{s=1}^{\varpi}\{P_s+y_s: y_s\in F\}\right|=x_F\varpi=x_F\prod_{j=1}^{i}(q^{j+e-1}+1),
$$
  Thus, ${\cal L}\cap{\cal O}_{\nu}(F)$ is  a Cameron-Liebler set in ${\cal O}_{\nu}(F)$ with parameter $x_F$.

  Let ${\cal L}\cap{\cal O}_{\nu}(F)$ be  a Cameron-Liebler set in ${\cal O}_{\nu}(F)$ with parameter $x_F$,
  and $P$ be a fixed maximal totally isotropic subspace contained in $F$, where $F=\langle e_1,e_2,\ldots,e_{\nu+i}\rangle$.
  Then $\{P+y : y\in F\}$ is a $(\nu,0)$-spread in $F$, and $\{P+y : y\in \mathbb{F}_q^{2\nu}\}$ is a $(\nu,0)$-spread in $\mathbb{F}_q^{2\nu}$.
  From ${\cal L}\cap{\cal O}_{\nu}(F)\cap\{P+y : y\in F\}={\cal L}\cap\{P+y : y\in F\}\subseteq{\cal L}\cap\{P+y : y\in \mathbb{F}_q^{2\nu}\}$,
  we deduce that $|{\cal L}\cap{\cal O}_{\nu}(F)\cap\{P+y : y\in F\}|=x_F\leq x=|{\cal L}\cap\{P+y : y\in \mathbb{F}_q^{2\nu}\}|$. Thus,
  $x_F\leq\min\{x,q^i\}$ since $x_F\leq q^i=|{\cal O}_{\nu}(F)|/\prod_{j=1}^{i}(q^{j+e-1}+1)$.
  $\qed$

The following result is obtained by using a method in \cite{De Beule2}.

\begin{thm}\label{lemma5.8}
Suppose that ${\cal L}$ is a Cameron-Liebler set in ${\cal O}_{\nu}$ with parameter $x$ and $S\in{\cal L}$.
Then $$x=\frac{\sum\limits_{S\subseteq T}x_T}{{\nu-1\brack \nu-i}_q}-\frac{q^\nu-1}{q^i-1}+1,$$
where $x_T$ is the parameter of the Cameron-Liebler set ${\cal L}\cap{\cal O}_\nu(T)$ in ${\cal O}_{\nu}(T)$ and the sum runs over all $(\nu+i,2i)$-flats containing $S$ in $ACG(2\nu,\,\mathbb{F}_q)$.
\end{thm}
\proof Let ${\cal F}$ be the set of all $S'\in{\cal L}$ satisfying $(S,S')\in R_{(0,1)}$, and
${\cal T}$ be the set of all $(\nu+i,2i)$-flats in $ACG(2\nu,\,\mathbb{F}_q)$ containing $S$.
Count the number $N$ of pairs $(S',T)$, where $S'\in{\cal F},T\in{\cal T}$ and $S'\subseteq T$.
For a fixed $T\in{\cal T}$,  by Theorem~\ref{thm5.6} and $S\in{\cal L}\cap{\cal O}_{\nu}(T)$, ${\cal L}\cap{\cal O}_{\nu}(T)$ is a Cameron-Liebler set in ${\cal O}_{\nu}(T)$
with parameter $x_T\geq1$, and the number of elements $S'\in{\cal L}\cap{\cal O}_{\nu}(T)$ satisfying $(S,S')\in R_{(0,1)}$ is $x_T-1$.
By \cite[Theorems~3.38,~5.37 and~6.43]{Wan}, $|{\cal T}|={\nu\brack \nu-i}_q$. Therefore, we have
$$N=\sum_{S\subseteq T}(x_T-1)=\sum_{S\subseteq T}x_T-{\nu\brack \nu-i}_q.$$
For a fixed $S'\in{\cal F}$, by Lemma~\ref{lemma2.1}, $S\vee S'$ is a $(\nu+1,2)$-flat. By \cite[Theorems~3.38,~5.37 and~6.43]{Wan} again,
both $S$ and $S'$
are contained in exactly ${\nu-1\brack \nu-i}_q$ many $(\nu+i,2i)$-flats. By Lemma~\ref{lemma3.10}, we have
$$|{\cal F}|=\frac{N}{{\nu-1\brack \nu-i}_q}=\frac{\sum_{S\subseteq T}x_T}{{\nu-1\brack \nu-i}_q}-\frac{{\nu\brack \nu-i}_q}{{\nu-1\brack \nu-i}_q}=x-1,$$
which implies that the desired result follows. $\qed$

\begin{cor}\label{lemma5.9}
Let ${\cal L}$ be a Cameron-Liebler set in ${\cal O}_{\nu}$ with parameter $x\geq2$ and $S\in{\cal L}$. Suppose that
${\cal T}$ is the set of all $(\nu+i,2i)$-flats in $ACG(2\nu,\,\mathbb{F}_q)$ containing $S$ and ${\cal T}_\theta=\{T\in{\cal T} : x_T=\theta\}$ for all $\theta=1,\ldots,\min\{x,q^i\}$,
where $x_T$ is the parameter of the Cameron-Liebler set ${\cal L}\cap{\cal O}_\nu(T)$ in ${\cal O}_{\nu}(T)$. Then
$$
{\nu\brack \nu-i}_q=\sum\limits_{\theta=1}^{\min\{x,q^i\}}|{\cal T}_\theta| \quad\hbox{and}\quad
(x-1){\nu-1\brack \nu-i}_q=\sum\limits_{\theta=2}^{\min\{x,q^i\}}(\theta-1)|{\cal T}_\theta|.
$$
\end{cor}
\proof
By Theorem~\ref{lemma5.8}, we have
$$x=\frac{\sum\limits_{T\in{\cal T}}(x_T-1)}{{\nu-1\brack \nu-i}_q}+1.$$
Since $F\in{\cal L}\cap{\cal O}_\nu(T)$, we obtain $x_T\geq 1$.
Since $|{\cal T}|={\nu\brack \nu-i}_q$ and $x_T\leq \min\{x,q^i\}$ for each $T\in{\cal T}$, we have
\begin{eqnarray*}
{\nu\brack \nu-i}_q&=&|{\cal T}|=\sum\limits_{\theta=1}^{\min\{x,q^i\}}|{\cal T}_\theta|, \\
(x-1){\nu-1\brack \nu-i}_q&=&\sum\limits_{T\in{\cal T}}(x_T-1)=\sum\limits_{\theta=2}^{\min\{x,q^i\}}(\theta-1)|{\cal T}_\theta|.
\end{eqnarray*}
Therefore, the desired result follows. $\qed$

\begin{cor}\label{thm5.10}
Let ${\cal L}$ be a Cameron-Liebler set in ${\cal O}_{\nu}$ with parameter $x\geq2$ and $S\in{\cal L}$.
Suppose that ${\cal T}$ and ${\cal T}_\theta \;(\theta=1,\ldots,\min\{x,q^i\})$ are as in Corollary~\ref{lemma5.9}. Then the following hold.
\begin{itemize}
\item[\rm(i)]
$|{\cal T}_1|\leq {\nu\brack\nu-i}_q-\frac{x-1}{\min\{x,q^i\}-1}{\nu-1\brack\nu-i}_q$.

\item[\rm(ii)]
If $|{\cal T}_1|={\nu\brack\nu-i}_q-\frac{x-1}{\min\{x,q^i\}-1}{\nu-1\brack\nu-i}_q$, then $|{\cal T}_{\min\{x,q^i\}}|=\frac{x-1}{\min\{x,q^i\}-1}{\nu-1\brack\nu-i}_q$
and $|{\cal T}_\theta|=0$ for all $\theta=2,\ldots,\min\{x,q^i\}-1$.

\item[\rm(iii)]
If $|{\cal T}_1|< {\nu\brack\nu-i}_q-\frac{x-1}{\min\{x,q^i\}-1}{\nu-1\brack\nu-i}_q$, then there exists some $1\leq\ell<\frac{q^\nu-1}{q^i-1}-\frac{x-1}{\min\{x,q^i\}-1}$ such that $(\ell-1){\nu-1\brack\nu-i}_q\leq|{\cal T}_1|<\ell{\nu-1\brack\nu-i}_q$ and $x\geq\frac{q^{\nu}-1}{q^i-1}-\ell+2$.
\end{itemize}
\end{cor}
\proof (i) Since $x_T\leq \min\{x,q^i\}$ for each $T\in{\cal T}$, by Corollary~\ref{lemma5.9}, we obtain
$$(\min\{x,q^i\}-1)\sum\limits_{\theta=2}^{\min\{x,q^i\}}|{\cal T}_\theta|\geq\sum\limits_{\theta=2}^{\min\{x,q^i\}}(\theta-1)|{\cal T}_\theta|=(x-1){\nu-1\brack\nu-i}_q,$$
which implies that $\sum\limits_{\theta=2}^{\min\{x,q^i\}}|{\cal T}_\theta|\geq\frac{x-1}{\min\{x,q^i\}-1}{\nu-1\brack\nu-i}_q$. By Corollary~\ref{lemma5.9} again, we obtain
$${\nu\brack \nu-i}_q=\sum\limits_{\theta=1}^{\min\{x,q^i\}}|{\cal T}_\theta|\geq|{\cal T}_1|+\frac{x-1}{\min\{x,q^i\}-1}{\nu-1\brack\nu-i}_q,$$
which implies that $|{\cal T}_1|\leq{\nu\brack\nu-i}_q-\frac{x-1}{\min\{x,q^i\}-1}{\nu-1\brack\nu-i}_q$.

(ii) If $|{\cal T}_1|={\nu\brack\nu-i}_q-\frac{x-1}{\min\{x,q^i\}-1}{\nu-1\brack\nu-i}_q$, then $\sum\limits_{\theta=2}^{\min\{x,q^i\}}|{\cal T}_\theta|=\frac{x-1}{\min\{x,q^i\}-1}{\nu-1\brack\nu-i}_q$, which implies that
$$(x-1){\nu-1\brack\nu-i}_q=(\min\{x,q^i\}-1)\sum\limits_{\theta=2}^{\min\{x,q^i\}}|{\cal T}_\theta|=\sum\limits_{\theta=2}^{\min\{x,q^i\}}(\theta-1)|{\cal T}_\theta|.$$
So, $|{\cal T}_{\min\{x,q^i\}}|=\frac{x-1}{\min\{x,q^i\}-1}{\nu-1\brack\nu-i}_q$
and $|{\cal T}_\theta|=0$ for all $\theta=2,\ldots,\min\{x,q^i\}-1$.

(iii) Since $|{\cal T}_1|< {\nu\brack\nu-i}_q-\frac{x-1}{\min\{x,q^i\}-1}{\nu-1\brack\nu-i}_q$ and ${\nu\brack \nu-i}/{\nu-1\brack\nu-i}_q=\frac{q^\nu-1}{q^i-1}$, there exists some $1\leq\ell<\frac{q^\nu-1}{q^i-1}-\frac{x-1}{\min\{x,q^i\}-1}$ such that $(\ell-1){\nu-1\brack\nu-i}_q\leq|{\cal T}_1|<\ell{\nu-1\brack\nu-i}_q$. By Corollary~\ref{lemma5.9}, we obtain
$${\nu\brack \nu-i}_q-|{\cal T}_1|=\sum\limits_{\theta=2}^{\min\{x,q^i\}}|{\cal T}_\theta|
\leq\sum\limits_{\theta=2}^{\min\{x,q^i\}}(\theta-1)|{\cal T}_\theta|=(x-1){\nu-1\brack\nu-i}_q,$$
which implies that
$$x\geq\frac{{\nu\brack \nu-i}_q}{{\nu-1\brack\nu-i}_q}-\frac{|{\cal T}_1|}{{\nu-1\brack\nu-i}_q}+1=\frac{q^{\nu}-1}{q^i-1}-\ell-\frac{|{\cal T}_1|-(\ell-1){\nu-1\brack\nu-i}_q}{{\nu-1\brack\nu-i}_q}+2.$$
It follows that $x\geq\frac{q^{\nu}-1}{q^i-1}-\ell+2$.
$\qed$

Note that ${\cal O}_{\nu}$ is a Cameron-Liebler set in ${\cal O}_{\nu}$ with parameter $x=q^\nu$.
Suppose that ${\cal T}$ and  ${\cal T}_\theta \;(\theta=1,\ldots,q^i)$ are as in Corollary~\ref{lemma5.9}. By Corollary~\ref{thm5.10} (ii),
$|{\cal T}_\theta|=0$ for all $\theta=1,\ldots,q^i-1$, and $|{\cal T}_{q^i}|=\frac{q^\nu-1}{q^i-1}{\nu-1\brack \nu-i}_q={\nu\brack \nu-i}_q$.
Since $|{\cal T}|={\nu\brack \nu-i}_q$,  ${\cal O}_{\nu}(T)$ is a Cameron-Liebler set in ${\cal O}_{\nu}(T)$ with parameter $x_T=q^i$ for every $T\in{\cal T}$.

\section{Concluding remarks}
Let $\mathbb{F}_q^{2\nu+\delta}$ be the $(2\nu+\delta)$-dimensional classical
space with parameter $e$,  and $ACG(2\nu+\delta,\mathbb{F}_q)$ be the corresponding $(2\nu+\delta)$-dimensional classical affine
space over a $q$-element finite field $\mathbb{F}_q$, where $(\delta,e)=(0,1)$ for the symplectic case, $(\delta,e)=(0,1/2)$ or $(1,3/2)$ for the unitary case, and $(\delta,e)=(0,0),(1,1)$ or $(2,2)$ for the orthogonal case, see \cite{Wan}.
Let ${\cal O}_{\nu}$ be the set of all maximal totally isotropic flats in $ACG(2\nu+\delta,\mathbb{F}_q)$.

Let $M$ be the incidence matrix with rows indexed with vectors in $\mathbb{F}_q^{2\nu+\delta}$
and columns indexed with flats in ${\cal O}_{\nu}$ such that entry $M_{x,F}=1$ if
and only if $x\in F$.
 A subset ${\cal L}$ of ${\cal O}_{\nu}$ is called a {\it Cameron-Liebler set} with parameter $x=|{\cal L}|/\prod_{i=1}^{\nu}(q^{i+e-1}+1)$
in ${\cal O}_{\nu}$ if $\chi_{{\cal L}}\in {\rm Im}(M^{\top})$.
In this paper, we only study Cameron-Liebler sets in ${\cal O}_{\nu}$ for $\delta=0$.
It seems interesting to discuss Cameron-Liebler sets in ${\cal O}_{\nu}$ for $\delta>0$.

Let $G$ be a finite group and $X$ an orbit under the action of $G$. If $R_0, R_1,\ldots, R_d$ are the orbits of $G$ on $X\times X$,
then $(X, \{R_i\}_{0\leq i\leq d})$ is an association scheme, see \cite{Bannai}.
By using this result, we know that ${\cal O}_{\nu}$ has a structure of an association scheme ${\mathfrak X}({\cal O}_{\nu})$.
In \cite{Guo,Liu,Liu2}, the scheme ${\mathfrak X}({\cal O}_{\nu})$ is constructed and  its character table is determined for $\delta=0$.
It seems interesting to discuss  the scheme ${\mathfrak X}({\cal O}_{\nu})$ for $\delta>0$.

\section*{ACKNOWLEDGMENT}
The authors are indebted to the anonymous reviewers for their detailed reports and constructive suggestions.
This research is supported by National Natural Science Foundation of China (11971146),  Foundation of Langfang Normal University (LSLB201707, LSPY201819, LSPY201915), Scientific Research Innovation Team of Langfang Normal University, and the Key Programs of Scientific Research Foundation of Hebei Educational Committee (ZD2019056).

\section*{CONFLICT OF INTEREST STATEMENT}
The author declares no conflict of interest.

\section*{DATA AVAILABILITY STATEMENT}
Data sharing is not applicable to this article as no datasets were generated or analysed during
the current study.

\section*{ORCID}
Jun Guo: http://orcid.org/0000-0001-9296-5608

\end{CJK*}

\end{document}